

\input amstex
\documentstyle{amsppt}

\input label.def
\input degt.def


\input epsf
\def\picture#1{\epsffile{#1-bb.eps}}
\def\cpic#1{$\vcenter{\hbox{\picture{#1}}}$}

\def\ie{\emph{i.e\.}}
\def\eg{\emph{e.g\.}}
\def\cf{\emph{cf\.}}
\def\via{\emph{via}}
\def\etc{\emph{etc}}

{\catcode`\@11
\gdef\proclaimfont@{\sl}}

\newtheorem{Remark}\Remark\thm\endAmSdef
\conjecture\thm\endproclaim
\def\paragraph{\subsection{}}

\def\dash{\item"\hfill--\hfill"}
\def\Dashes{\widestnumber\item{--}\roster}
\def\endDashes{\endroster}

\loadbold
\def\bA{\bold A}
\def\bD{\bold D}
\def\bE{\bold E}
\def\bJ{\bold J}

\def\tA#1{\smash{\tilde\bA#1}}
\def\tD#1{\smash{\tilde\bD#1}}
\def\tE#1{\smash{\tilde\bE#1}}
\def\tJ#1{\smash{\tilde\bJ#1}}
\def\I{{\fam0I}}
\def\II{{\fam0II}}
\def\III{{\fam0III}}
\def\IV{{\fam0IV}}

\let\Gk\kappa

\let\splus\oplus

\def\F{\Bbb F}

\def\CG#1{\Z_{#1}}
\def\BG#1{\Bbb B_{#1}}

\def\DG#1{\Bbb D_{#1}}
\def\SGSet{\Bbb S}
\def\SG#1{\SGSet_{#1}}

\let\Gs\sigma

\def\SL{\operatorname{\text{\sl SL\/}}}

\def\PSL{\operatorname{\text{\sl PSL\/}}}

\def\CC{\Cal C}

\let\B=B
\def\BB{\topsmash{\tilde B}} 
\def\kk{\topsmash{\tilde k}}            

\def\bR{\bar R}

\let\CC=C         
\let\branch\Gamma
\def\bbranch{\topsmash{\tilde\branch}}

\def\Cp#1{\Bbb P^{#1}}
\def\term#1-{$\DG{#1}$-}

\let\Ga\alpha
\let\Gb\beta
\let\Gg\gamma
\let\Gd\delta

\let\Gs\sigma
\let\Gr\rho

\def\1{^{-1}}

\def\0{{\setbox0\hbox{0}\hbox to\wd0{\hss\rm--\hss}}}

\def\ls|#1|{\mathopen|#1\mathclose|}
\let\<\langle
\let\>\rangle
\let\onto\twoheadrightarrow

\def\Aut{\operatorname{Aut}}
\def\Conv{\operatorname{conv}}

\def\ord{\operatorname{ord}}
\def\mult{\operatorname{mult}}
\def\Sk{\operatorname{Sk}}
\def\SkB{\Sk_{\B}}

\let\MB=M

\def\tabstrut{\vrule height9.5pt depth2.5pt}
\def\exstrutii{\omit\vrule height2pt\hss}
\def\exstrut{\exstrutii\vrule}

\def\inserthyphen{\ifcat\next a-\fi\ignorespaces}
\let\BLACK\bullet
\let\WHITE\circ
\def\CROSS{\vcenter{\hbox{$\scriptstyle\mathord\times$}}}
\let\STAR*
\def\pblack-{$\BLACK$\futurelet\next\inserthyphen}
\def\pwhite-{$\WHITE$\futurelet\next\inserthyphen}
\def\pcross-{$\CROSS$\futurelet\next\inserthyphen}
\def\pstar-{$\STAR$\futurelet\next\inserthyphen}
\def\black{\protect\pblack}
\def\white{\protect\pwhite}

\def\star{\protect\pstar}
\def\NO#1{\mathord\#_{#1}}
\def\nblack{\NO\BLACK}
\def\nwhite{\NO\WHITE}

\def\nstar{\NO\STAR}

\def\No{no\.~}
\def\Nos{nos\.~}

\def\GAP{{\tt GAP}}
\def\FRAG#1{--$#1$}
\def\fragment(#1){\ref{fig.d}(#1)\FRAG}
\def\frag{\ref{fig.d5}\FRAG}

\def\beginGAP{\bgroup
 \catcode`\^=12\catcode`\#=12\catcode`\_=12
 \obeylines\obeyspaces\eightpoint
 \tt}
\let\endGAP\egroup

\def\Tab{\vtop\bgroup\openup1pt\halign\bgroup$##$\hss\cr}
\def\endTab{\crcr\egroup\egroup}

\def\gm{\frak m}
\def\tgm{\tilde\gm}
\let\disk\Delta
\let\slope\varkappa
\def\ts{\bar\slope}

\def\tsp{\operatorname{ts}}
\def\tspB{\tsp_{\B}}
\def\tdf{\operatorname{td}}
\def\tdfB{\tdf_{\B}}

\def\ixi{_{\fam0I}}
\def\ixii{_{\fam0II}}
\def\Fi{F\ixi}
\def\Fii{F\ixii}
\def\Ri{R\ixi}
\def\Rii{R\ixii}
\def\vi{v\ixi}
\def\vii{v\ixii}

\def\corner#1#2{\,\widehat{\!#1#2\!}\,}

\topmatter

\author
Alex Degtyarev
\endauthor

\title
The fundamental group of a generalized trigonal curve
\endtitle

\address
Department of Mathematics,
Bilkent University,
06800 Ankara, Turkey
\endaddress

\email
degt\@fen.bilkent.edu.tr
\endemail

\abstract
We develop a modification of the Zariski--van Kampen approach for
the computation of the fundamental group of a trigonal curve with
improper fibers. As an application, we list the deformation
families and compute the fundamental groups of all irreducible
maximizing simple sextics with a type~$\bD$ singular point.
\endabstract

\keywords
Plane sextic, fundamental group, trigonal curve, dessin d'enfant
\endkeywords

\subjclassyear{2000}
\subjclass
Primary: 14H45; 
Secondary: 14H30, 
14H50 
\endsubjclass

\endtopmatter

\document

\section{Introduction}

\subsection{Principal results}
We attempt to develop a modification of the classical Zariski--van
Kampen approach~\cite{vanKampen} suitable to compute the
fundamental group of a generalized trigonal curve, \ie, a trigonal
curve with \emph{improper fibers}, at which the curve meets the
exceptional section. A similar question was addressed
in~\cite{dessin.e7}, where the only improper fiber was `hidden' at
infinity. Here, we consider the case of arbitrarily many improper
fibers (up to two in the applications).

The basic tool used in Zariski--van Kampen's method is the braid
monodromy related to an appropriate pencil. In order to keep the
braid monodromy well defined, $\BG3$-valued, and easily computable
\via\ skeletons (see Subsection~\ref{s.braid.Sk}), we pass to the
associated genuine trigonal curve and introduce the concept of
\emph{slopes}, which compensate for the improper fibers. We
compute local slopes (Subsection~\ref{s.local}), study their
properties,
and discuss the
modifications that should be made to the braid
relations~\eqref{rel.braid} and
relation at infinity~\eqref{rel.infinity}
in the Zariski--van Kampen presentation of the fundamental group,
see Corollary~\ref{vK.slopes}.

As a simple application, in
Subsection~\ref{s.quintics} we recompute the fundamental groups of
irreducible plane quintics with a double point. (These groups were
originally found in~\cite{groups} and~\cite{Artal}, but the
computation \via\ trigonal curves is much simpler and more
straightforward; it could easily be computerized.)

\subsection{Plane sextics}\label{s.secxtics}
A more advanced example is the case of
irreducible plane sextics with a
type~$\bD$ singular point.

Recall that a plane sextic $\CC\subset\Cp2$ is called
\emph{simple} if all its singular points are simple, \ie, those of
types~$\bA_p$, $\bD_q$, $\bE_6$, $\bE_7$, or~$\bE_8$ (see
\eg~\cite{AVG} for the notation). The total Milnor number $\mu(\CC)$
of a simple sextic~$C$ does not exceed~$19$; if $\mu(\CC)=19$, the
sextic is called \emph{maximizing}. Maximizing sextics are rigid:
if two such sextics are equisingular deformation equivalent, they
are related by a projective transformation. Each maximizing sextic
is defined over an algebraic number field.

A sextic is said to be of \emph{torus type} if its equation can be
represented in the form $f_2^3+f_3^2=0$, where $f_2$ and~$f_3$ are
certain homogeneous polynomials of degree~$2$ and~$3$,
respectively. Alternatively, $\CC$ is of torus type if it is the
ramification locus of a projection to~$\Cp2$ of a cubic surface
$V\subset\Cp3$. This property is invariant under equisingular
deformations. Each sextic~$\CC$ of torus type can be perturbed to a six
cuspidal sextic, see~\cite{Zariski.group}, hence the
fundamental group
$\pi_1(\Cp2\sminus\CC)$ factors to the reduced braid group
$\bar\BG3:=\BG3/(\Gs_1\Gs_2)^3\cong\CG2*\CG3\cong\PSL(2,\Z)$; in
particular, this group is never abelian or finite.

In this paper, we study irreducible maximizing simple sextics with
a type~$\bD$ singular point and without type~$\bE$ singular
points. (Sextics with type~$\bE$ points are the subject
of~\cite{dessin.e7}, \cite{dessin.e8}, and~\cite{dessin.e6}.) We
list the equisingular deformation families of such sextics
(Theorem~\ref{th.d})
and compute their fundamental groups (Theorem~\ref{th.group}).
As in the previous papers, the
principal tool is the reduction of a sextic with a triple singular
point to a generalized trigonal curve in~$\Sigma_1$.

\theorem\label{th.d}
There are $38$ deformation families of irreducible maximizing
simple sextics with a type~$\bD$ singular point
and without type~$\bE$ singular points, realizing $25$
sets of singularities \rom(see Tables~\ref{tab.d}
and~\ref{tab.d5} in Section~\ref{S.proof.d}\rom).
One of the families is of torus type
\rom(the set of singularities
$\bD_{5}\splus(\bA_8\splus3\bA_2)$, \No$27$ in
Table~\ref{tab.d5}\rom)\rom;
the others are not.
\endtheorem

Theorem~\ref{th.d} is proved in Section~\ref{S.proof.d}. In
principle, the statement can be obtained by comparing the results
of J.-G.~Yang~\cite{Yang} (a list of all sets of singularities
that can be
realized by
an irreducible maximizing simple sextic) and
I.~Shimada~\cite{Shimada}
(a list of sets of singularities represented by several
deformation families), using the global Torelli theorem for
$K3$-surfaces. The advantage of our approach is an explicit
construction of each sextic, which can further be used in the
study of its geometry.

\theorem\label{th.group}
Let~$\CC\subset\Cp2$ be
an irreducible maximizing simple sextic with a
type~$\bD$ singular point.
If $\CC$ is of torus type, then
$\pi_1(\Cp2\sminus \CC)$ is the reduced braid group
$\bar\BG3=\BG3/(\Gs_1\Gs_2)^3\cong\CG2*\CG3$\rom; otherwise,
$\pi_1(\Cp2\sminus \CC)=\CG6$.
\endtheorem

If $\CC$ has a type~$\bE$ point, the statement follows
from~\cite{dessin.e7},
\cite{dessin.e8}, and~\cite{dessin.e6}. Other sextics as in
Theorem~\ref{th.group} are considered in
Section~\ref{S.proof.group}, using the models constructed in
Section~\ref{S.proof.d} and the
approach developed in
Section~\ref{S.slopes}.
As an immediate consequence, one obtains the following corollary.

\corollary
Let $\CC'$ be a perturbation of a sextic~$\CC$ as in
Theorem~\ref{th.group}. If $\CC'$ is of torus type, then
$\pi_1(\Cp2\sminus \CC')=\bar\BG3$\rom; otherwise,
$\pi_1(\Cp2\sminus \CC')=\CG6$.
\qed
\endcorollary

Recall that any induced subgraph of the combined Dynkin graph of a
simple sextic~$C$ can be realized by a perturbation of~$C$.

We do not treat systematically reducible curves, as that would
require an enormous amount of work. However, as a simple
by-product, we do compute the groups of a few
maximizing deformation
families and their perturbations,
see Table~\ref{tab.d-r} in Subsection~\ref{s.RR0}
and Table~\ref{tab.d5-r} in Subsection~\ref{s.group.D5}.
Perturbing, one obtains more irreducible sextics with abelian
groups, see Proposition~\ref{perturbations}.
Altogether, the results of this and a few previous papers suggest
the following conjecture.

\conjecture
With the exception of the maximizing sextics realizing
the following three sets of singularities\rom:
\Dashes
\dash
$2\bE_6\splus\bA_4\splus\bA_3$ \rom(two curves\rom;
$\pi_1=\SL(2,\F_5)\rtimes\CG6$, see~\cite{dessin.e6}\rom),
\dash
$\bE_7\splus2\bA_4\splus2\bA_2$ \rom(one curve\rom;
$\pi_1=\SL(2,\F_{19})\rtimes\CG6$, see~\cite{dessin.e7}\rom),
\dash
$\bE_8\splus\bA_4\splus\bA_3\splus2\bA_2$ \rom(one curve\rom;
$\pi_1=\SL(2,\F_5)\odot\CG{12}$, see~\cite{dessin.e8}\rom),
\endDashes
the fundamental group $\pi_1:=\pi_1(\Cp2\sminus\CC)$ of an
irreducible simple
sextic $\CC\subset\Cp2$ that is not of torus type and has a
triple singular point is abelian.
\endconjecture

(In the description of the groups, $\rtimes$ stands for a
semi-direct product and $\odot$ stands for a central product:
$\SL(2,\F_5)\odot\CG{12}$ is the quotient of
$\SL(2,\F_5)\times\CG{12}$ by the diagonal subgroup
$\CG2\subset\operatorname{Center}\SL(2,\F_5)\times\CG2$.)
A proof of this conjecture would require a detailed study of the
degenerations, which would probably lead to reducible sextics, and
a computation of the groups of (some) reducible maximizing sextics
with a type~$\bD$ or type~$\bE_7$ singular point.

After Theorem~\ref{th.group},
there still remain five maximizing simple sextics of torus type
with unknown fundamental groups;
their sets of singularities are
$$
\gathered
(\bA_{14}\splus\bA_2)\splus\bA_3,\quad
(\bA_{14}\splus\bA_2)\splus\bA_2\splus\bA_1,\quad
(\bA_{11}\splus2\bA_2)\splus\bA_4,\\
(\bA_8\splus\bA_5\splus\bA_2)\splus\bA_4,\quad
(\bA_8\splus3\bA_2)\splus\bA_4\splus\bA_1.
\endgathered
$$
(We use the list of
irreducible sextics of torus type found in~\cite{OkaPho.moduli};
maximizing sets of singularities can also be extracted
from~\cite{Yang}.
Due to~\cite{Shimada},
$(\bA_8\splus\bA_5\splus\bA_2)\splus\bA_4$ is realized by a pair
of complex conjugate curves, whereas the four
remaining sets of singularities define a single
deformation family each.)
Assuming that, up to complex conjugation, each
\emph{non}-maximizing set of singularities is realized by at most
one connected deformation family of sextics of torus type (which
is probably true, but proof is still pending), the groups of all
such sextics are known.
For details and further references, see recent
survey~\cite{survey}.

\subsection{Contents of the paper}
In Section~\ref{S.trigonal}, we introduce the terminology and
remind a few known results related to generalized trigonal curves.
Section~\ref{S.slopes} deals with the fundamental groups: we
remind the general approach, due to Zariski and van
Kampen~\cite{vanKampen}, specialize it to genuine trigonal curves
(following~\cite{degt.e6}),
and introduce slopes for generalized trigonal curves. Then, we
explain how the slopes and the global monodromy can be computed
and consider an example, applying the approach to irreducible
plane quintics. In Section~\ref{S.proof.d}, we enumerate the
deformation families of sextics as in Theorem~\ref{th.d} by
describing the skeletons of their trigonal models; this
description is used in Section~\ref{S.proof.group} in the
computation of the fundamental groups.

\subsection{Acknowledgements}
I am grateful to E.~Artal Bartolo, who helped me to identify the
group~$\bar\BG3$
of the sextic of torus type in Theorem~\ref{th.group}.

\section{Generalized trigonal curves\label{S.trigonal}}

In this section, we mainly introduce the terminology and cite a
few known results related to (generalized) trigonal curves in
Hirzebruch surfaces. Principal references
are~\cite{degt.kplets} and~\cite{dessin.e7}.

\subsection{Hirzebruch surfaces}\label{s.surface}
Recall that the \emph{Hirzebruch surface} $\Sigma_k$, $k\ge0$, is
a rational geometrically ruled surface with an
\emph{exceptional section} $E=E_k$ of
self-intersection~$-k$. The fibers
of the ruling are referred to as the fibers of~$\Sigma_k$.
The semigroup of classes of effective divisors
on
~$\Sigma_k$
is generated by the classes of the exceptional
section~$E$ and a fiber~$F$; one has
$E^2=-k$, $F^2=0$, and $E\cdot F=1$.

Fix a Hirzebruch surface~$\Sigma_k$, $k\ge1$. Denote by
$p\:\Sigma_k\to\Cp1$ the ruling, and let $E\subset\Sigma_k$ be
the exceptional section, $E^2=-k$. Given a point~$b$ in the
base~$\Cp1$, we denote by~$F_b$ the fiber $p^{-1}(b)$. (With a
certain abuse of the language, the points in the base~$\Cp1$ of
the ruling are also referred to as fibers of~$\Sigma_k$.)
Let
$F_b^\circ$ be the `open fiber' $F_b\sminus E$. Observe that
$F_b^\circ$ is a dimension~$1$ affine space over~$\C$; hence, one
can speak about lines, circles, angles, convexity, \etc\.
in~$F_b^\circ$.
In particular, one can define the \emph{convex hull} $\Conv S$
of a subset $S\subset\Sigma_k\sminus E$
as the union of its fiberwise convex
hulls:
$$\tsize
\Conv S=\bigcup_{b\in\Cp1}\Conv(S\cap F_b^\circ).
$$

\subsection{Trigonal curves}\label{s.trigonal}
A \emph{generalized trigonal curve} on a
Hirzebruch surface~$\Sigma_k$ is a reduced curve~$\B$
not containing the exceptional section~$E$ and
intersecting each generic fiber at three
points. In this paper, we assume in addition
that a trigonal curve does not contain a fiber of~$\Sigma_k$
as a component.

A \emph{singular fiber} of a generalized trigonal curve
$\B\subset\Sigma_k$ is a fiber~$F$ of~$\Sigma_k$ that is not
transversal to the union $\B\cup E$. Thus, $F$ is either
the fiber over a critical value of the
restriction to~$\B$ of the ruling $\Sigma_k\to\Cp1$
or
the fiber through a point of intersection
of~$\B$ and~$E$. In the former case, the fiber is called
\emph{proper}; in the latter case, the fiber is called
\emph{improper} and the points of intersection of~$\B$ and~$E$
are called \emph{points at infinity}. In general, the local
branches of~$\B$ that intersect a fiber~$F$ outside of~$E$ are
called \emph{proper} at~$F$.

A \emph{\rom(genuine\rom) trigonal curve}
is a generalized trigonal curve $\B\subset\Sigma_k$
disjoint from
the exceptional section.
One has $\B\in\ls|3E+3kF|$; conversely, any reduced curve
$\B\in\ls|3E+3kF|$ not containing~$E$ as a component is a trigonal
curve.

We use the following notation for the topological types of
proper fibers:
\Dashes
\dash
$\tA{_0}$: a nonsingular fiber;
\dash
$\tA{_0^*}$: a simple vertical tangent;
\dash
$\tA{_0^{**}}$: a vertical inflection tangent;
\dash
$\tA{_1^*}$: a node of~$\B$ with one of the branches vertical;
\dash
$\tA{_2^*}$: a cusp of~$\B$ with vertical tangent;
\dash
$\tA{_p}$, $p\ge2$, $\tD{_q}$, $q\ge4$,
$\tE{_{6r+\epsilon}}$, $r\ge1$, $\epsilon=0,1,2$,
$\tJ{_{r,p}}$, $r\ge2$, $p\ge0$:
a singular point of~$\B$ of the same type (see~\cite{AVG} for the
notation)
with minimal possible local intersection index with the fiber.
\endDashes
For `simple' fibers of types~$\tA{}$, $\tD{}$, $\tE{_6}$,
$\tE{_7}$, and~$\tE{_8}$, this notation refers to the incidence
graph of $(-2)$-curves in the corresponding singular elliptic fiber;
this graph is an affine Dynkin diagram.

\Remark
The topological classification of singular fibers of trigonal
curves is close to that for elliptic surfaces, see~\cite{Kodaira},
except that in this paper we admit curves with non-simple
singularities. It would probably be more convenient (but slightly
less transparent) to use an appropriate extension of Kodaira's
notation, for example $\I_p^r$, $\II^r$, $\III^r$, and~$\IV^r$,
with $r=0$ and~$1$ referring, respectively, to the empty subscript
and~$^*$ in~\cite{Kodaira}. Another alternative would be to extend
the series $\tJ{_{r,p}}$ and $\tE{_{6r+\epsilon}}$ to the values
$r=0$ and~$1$. Among other advantages, in both cases an elementary
transformation (see Subsection~\ref{s.transform} below) would
merely increase the value of~$r$ by~$1$. However, I chose to
retain the commonly accepted notation for the types of simple
singularities.
\endRemark

The fibers of types~$\tA{_0^{**}}$, $\tA{_1^*}$, and~$\tA{_2^*}$
are called \emph{unstable}; all other singular fibers are called
\emph{stable}. A trigonal curve~$\B$ is \emph{stable} if so are all its
singular fibers. (This notion of stability differs from the one
accepted in algebraic geometry; we refer to the topological stability
under equisingular deformations of~$\B$. An unstable fiber may
split as follows: $\tA{_0^{**}}\to2\tA{_0^*}$,
$\tA{_1^*}\to\tA{_1}\splus\tA{_0^*}$, or
$\tA{_2^*}\to\tA{_2}\splus\tA{_0^*}$, the splitting not changing
the topology of the pair $(\Sigma_k,\B)$.)

The \emph{multiplicity} $\mult F$ of a singular fiber~$F$ of a
trigonal curve~$\B$ is the number of simplest (\ie,
type~$\tA{_0^*}$) fibers into which $F$ splits under deformations
of~$\B$. For the $\tA{}$ type fibers, one has
$\mult\tA{_0}=0$,
$\mult\tA{_0^*}=1$,
$\mult\tA{_0^{**}}=2$,
$\mult\tA{_1^*}=3$,
$\mult\tA{_2^*}=4$,
and $\mult\tA{_{p}}=p+1$ for $p>0$. Each elementary transformation
(see Subsection~\ref{s.transform} below) contracting~$F$ increases
$\mult F$ by~$6$. The sum of the multiplicities of all singular
fibers of a trigonal curve $\B\subset\Sigma_k$ equals $12k$.

\subsection{Elementary transformations}\label{s.transform}
An \emph{elementary transformation} of~$\Sigma_k$
is a birational transformation
$\Sigma_k\dashrightarrow\Sigma_{k+1}$ consisting in blowing up a
point~$P$ in the exceptional section of~$\Sigma_k$ followed by blowing
down the fiber~$F$ through~$P$. The inverse transformation
$\Sigma_{k+1}\dashrightarrow\Sigma_{k}$
blows up
a point~$P'$
\emph{not} in the exceptional section of~$\Sigma_{k+1}$
and blows down
the fiber~$F'$ through~$P'$.

An elementary transformation converts a proper fiber as follows:
\roster
\item\local{A0}
$\tA{_0}\to\tD{_4}\to\tJ{_{2,0}}\to\ldots\to\tJ{_{r,0}}\to\ldots$
\ (not detected by the $j$-invariant);
\item\local{A0*}
$\tA{_0^*}\to\tD{_5}\to\tJ{_{2,1}}\to\ldots\to\tJ{_{r,1}}\to\ldots$
\ ($j=\infty$, $\ord j=1$);
\item\local{Ap}
$\tA{_{p-1}}\to\tD{_{p+4}}\to\tJ{_{2,p}}\to\ldots\to\tJ{_{r,p}}\to\ldots$
\ ($p\ge2$; $j=\infty$, $\ord j=p$);
\item\local{A0**}
$\tA{_0^{**}}\to\tE{_6}\to\tE{_{12}}\to\ldots\to\tE{_{6r}}\to\ldots$
\ ($j=0$, $\ord j=1\bmod3$);
\item\local{A1*}
$\tA{_1^*}\to\tE{_7}\to\tE{_{13}}\to\ldots\to\tE{_{6r+1}}\to\ldots$
\ ($j=1$, $\ord j=1\bmod2$);
\item\local{A2*}
$\tA{_2^*}\to\tE{_8}\to\tE{_{14}}\to\ldots\to\tE{_{6r+2}}\to\ldots$
\ ($j=0$, $\ord j=2\bmod3$).
\endroster
For the reader's convenience, we also indicate the value $j=v$
and the ramification index $\ord j$ of the $j$-invariant, see
Subsection~\ref{s.j} below, which is invariant under elementary
transformations. In a neighborhood of the fiber, the $j$-invariant
has the form $v+t^{\ord j}$ if $v=0$ or~$1$ or $1/t^{\ord j}$ if
$ v=\infty$.

Let $\BB\subset\Sigma_{\kk}$ be a generalized trigonal curve.
Then, by a sequence of elementary transformations, one can resolve
the points of intersection of~$\BB$ and~$E$ and obtain a genuine
trigonal curve $\B\subset\Sigma_{k}$, $k\ge\kk$, birationally
equivalent to~$\BB$. The trigonal curve~$\B$ obtained from~$\BB$
by a minimal number of elementary transformations is called the
\emph{trigonal model} of~$\BB$.

\Remark\label{rem.simplified}
Alternatively, given a trigonal curve
$\B\subset\Sigma_k$ with triple singular points, one can apply a
sequence of inverse elementary transformations to
obtain a trigonal curve
$\B'\subset\Sigma_{k'}$, $k'\le k$, birationally
equivalent to~$\B$ and with $\tA{}$~type singular fibers only.
This curve~$\B'$ is called in~\cite{degt.kplets} the
\emph{simplified model} of~$\B$.
\endRemark

\subsection{The $j$-invariant}\label{s.j}
The \emph{\rom(functional\rom) $j$-invariant}
$j_{\B}\:\Cp1\to\Cp1$ of
a generalized
trigonal curve $\B\subset\Sigma_2$ is defined as the analytic
continuation of the function sending a point
$b$ in the base $\Cp1$ of~$\Sigma_2$
representing a
nonsingular fiber~$F_b$ of~$\B$
to the $j$-invariant (divided
by~$12^3$)
of the elliptic
curve covering~$F_b$
and ramified at the four points of intersection of~$F_b$
and $\B+E$. The curve~$\B$
is called \emph{isotrivial} if $j_{\B}=\const$. Such curves can
easily be enumerated, see \eg~\cite{degt.kplets}.

By definition, $j_{\B}$ is invariant under elementary
transformations. The values of~$j_{\B}$ at the singular fibers
of~$\B$ are listed in Subsection~\ref{s.transform}. The points
$b\subset\Cp1$ with $j_{\B}(b)=0$ and $\ord_b j_{\B}=0\bmod3$
or $j_{\B}(b)=1$ and $\ord_b j_{\B}=0\bmod2$ correspond to
fibers~$F_b$ admitting extra symmetries. Assuming~$F_b$ proper
(hence nonsingular), consider the three points of
intersection of~$\B$ and~$F^\circ_b$. Then
\Dashes
\dash
the three points form an equilateral triangle if
$j_{\B}(b)=0$, $\ord_b j_{\B}=0\bmod3$;
\dash
one of the points is at the center of the segment connecting
the two others if
$j_{\B}(b)=1$, $\ord_b j_{\B}=0\bmod2$.
\endDashes

\definition\label{def.max}
A non-isotrivial trigonal curve~$\B$
is called \emph{maximal} if it has the following properties:
\roster
\item\local{noD4}
$\B$ has no singular fibers of type~$\tD{_4}$ or
$\tJ{_{r,0}}$, $r\ge2$;
\item\local{0,1,infty}
$j=j_{\B}$ has no critical values other than~$0$, $1$, and~$\infty$;
\item\local{le3}
each point in the pull-back $j^{-1}(0)$ has ramification index at
most~$3$;
\item\local{le2}
each point in the pull-back $j^{-1}(1)$ has ramification index at
most~$2$.
\endroster
\enddefinition

An important property of maximal trigonal
curves is their rigidity, see~\cite{degt.kplets}: any small
fiberwise equisingular deformation of such a curve
$\B\subset\Sigma_k$ is
isomorphic to~$\B$. Any maximal trigonal curve is defined over an
algebraic number field. Such curves are classified by their
skeletons, see Theorem~\ref{th.skeleton} below.

A maximal trigonal curve~$\B$ with simple singularities only can
be characterized in terms of its total Milnor number $\mu(\B)$
(\ie, the sum of the Milnor numbers of all singular points
of~$\B$). The following criterion is proved in~\cite{dessin.e7}.

\theorem\label{th.maximal}
For a non-isotrivial genuine trigonal curve $\B\subset\Sigma_k$
with simple singularities only one has
$$
\mu(\B)\le5k-2-\#\{\text{\rm unstable fibers of~$\B$}\},
\eqtag\label{eq.maximal}
$$
the equality holding if and only if $\B$ is maximal.
\qed
\endtheorem

\Remark
The inequality in Theorem~\ref{th.maximal}
may not hold is $\B$ has non-simple singular
points, as each elementary transformation producing a non-simple
singular point increases~$\mu$ by~$6$ while increasing~$k$
by~$1$.
\endRemark

\subsection{Skeletons}\label{s.skeleton}
The \emph{skeleton} $\Sk=\SkB$ of a trigonal curve
$\B\subset\Sigma_k$ is defined as Grothendieck's \emph{dessin
d'enfants} of its $j$-invariant~$j_{\B}$. More precisely, $\Sk$ is
the planar map $j_{\B}\1([0,1])\subset S^2\cong\Cp1$. The
pull-backs of~$0$ are called \emph{\black-vertices}, and the
pull-backs of~$1$ are called \emph{\white-vertices}. The \black--
and \white-vertices are called \emph{essential}; the other
vertices that $\Sk$ may have (due to the critical values
of~$j_{\B}$ in the interval $(0,1)$\,) are called
\emph{unessential}.

By definition, $\Sk$ is a graph in the base of the ruling
$\Sigma_k\to\Cp1$, so that one can speak about the fibers
of~$\Sigma_k$ represented by points of~$\Sk$. On the other hand,
for the classification statements, see \eg\
Theorem~\ref{th.skeleton} below, it is important that $\Sk$ is
regarded as a graph in the \emph{topological} sphere~$S^2$; the
analytic structure is given by the skeleton itself \via\ Riemann's
existence theorem.

The \black-vertices
of valency $1\bmod3$ or $2\bmod3$ and
\white-vertices of valency $1\bmod2$ are called \emph{singular};
they correspond to singular fibers of the curve of one of the
types \iref{s.transform}{A0**}--\ditto{A2*}. All other \black--
and \white-vertices are called \emph{nonsingular}.

After a small fiberwise equisingular deformation of a trigonal
curve~$\B$ one can assume that its skeleton~$\SkB$ has the
following properties:
\roster
\item\local{essential}
all vertices of~$\SkB$ are essential;
\item\local{<3}
each \black-vertex has valency at most~$3$;
\item\local{<2}
each \white-vertex has valency at most~$2$.
\endroster
A skeleton satisfying these conditions is called \emph{generic}.
Note that any skeleton satisfying condition~\loccit{essential} is
a bipartite graph. For this reason, in the drawings below we omit
bivalent \white-vertices, assuming that such a vertex is to be
inserted in the middle of each edge connecting two
\black-vertices. In particular, for a generic skeleton, only
singular monovalent \white-vertices are drawn.

A \emph{region} of a skeleton $\Sk\subset\Cp1$ is
a connected component of the complement $\Cp1\sminus\Sk$.
One can also speak about \emph{closed regions}, which are
connected components of the manifold theoretical cut of~$\Cp1$
along~$\Sk$. (In general, a closed region $\bR$ is \emph{not} the
same as the closure of the corresponding open region~$R$.) We say
that a region~$R$ is an \emph{$m$-gon} (or an \emph{$m$-gonal
region})
if the boundary of the
corresponding closed region~$\bR$ contains $m$ \black-vertices.
For example, in Figure~\ref{fig.d}(b) below, the three regions
marked with~$\Ga$, $\Gb$, and $\bar\Gb$ are monogons, whereas the
outer region is a nonagon. In Figure~\ref{fig.d}(c), there are two
monogons (marked with~$\Ga$ and~$\Gb$) and two pentagons.

Each region~$R$ of~$\SkB$ contains a finite number of singular
fibers of~$\B$, which can be of one of the types
\iref{s.transform}{A0}--\ditto{Ap} (excluding $\tA{_0}$, which is
not singular). One can use a sequence of inverse elementary
transformations and convert these fibers to the $\tA{}$ type
fibers starting the series. If $R$ is an $m$-gonal region, the
total multiplicity of these $\tA{}$ type fibers equals~$m$.

\subsection{Skeletons and maximal curves}\label{s.maximal}
The skeleton~$\SkB$ of a maximal trigonal curve
$\B\subset\Sigma_k$ is
necessarily generic and connected. (It follows that each region
of~$\SkB$ is a topological disk.) Each $m$-gonal region~$R$
of~$\SkB$ contains a single singular fiber~$F_R$ of~$\B$; its type
is one of~\iref{s.transform}{A0*} if $m=1$ or one
of~\iref{s.transform}{Ap} with $p=m$ if $m\ge2$. Thus, the type
of~$F_R$ is determined by its multiplicity. The other singular
fibers of~$\B$ are over the singular vertices of~$\SkB$; the type
of such a singular fiber~$F_v$
is also determined by its
multiplicity (and the type and the valency of~$v$).

The function~$\tspB$ sending each region~$R$ to the
multiplicity $\mult F_R$ and each singular vertex~$v$ to the
multiplicity $\mult F_v$ is called the \emph{type
specification}. It has the following properties:
\roster
\item\local{ts.R}
$\tspB(\text{$m$-gonal region~$R$})=m+6s$, $s\in\Z_{\ge0}$;
\item\local{ts.black}
$\tspB(\text{singular \black-vertex~$v$})=2(\text{valency of~$v$})+6s$,
$s\in\Z_{\ge0}$;
\item\local{ts.white}
$\tspB(\text{singular \white-vertex})=3+6s$, $s\in\Z_{\ge0}$;
\item\local{ts.sum}
the sum of all values of~$\tspB$ equals $12k$.
\endroster

The following statement is essentially contained
in~\cite{degt.kplets}.

\theorem\label{th.skeleton}
The map $\B\mapsto(\SkB,\tspB)$
establishes a bijection between the set of isomorphism classes
\rom(equivalently, fiberwise equisingular deformation classes\rom)
of maximal trigonal curves in~$\Sigma_k$ and the set of
orientation preserving diffeomorphism classes of pairs
$(\Sk,\tsp)$, where $\Sk\subset S^2$ is a connected generic
skeleton and $\tsp$ is a function on the set of regions and
singular vertices of~$\Sk$
satisfying conditions
\itemref{s.maximal}{ts.R}--\ditto{ts.sum} above.
\qed
\endtheorem

\Remark\label{rem.td}
Often it is more convenient to replace~$\tspB$ with the
$\Z_{\ge0}$-valued function $\tdfB$ sending each region and
singular vertex to the integer~$s$ appearing
in~\itemref{s.maximal}{ts.R}--\ditto{ts.white}. In term
of~$\tspB$,
the index~$k$ of the Hirzebruch surface~$\Sigma_k$
containing~$\B$ is given as follows, \cf~\cite{dessin.e7}:
$$\tsize
\nblack+\nwhite(1)+\nblack(2)=2\bigl(k-\sum\tdfB\bigr),
$$
where $\nblack$ is the total number of \black-vertices,
$\nstar(i)$ is the number of \star-vertices of valency~$i$, and
$\sum\tdfB$ is the sum of all values of~$\tdfB$. The singular
points of~$\B$ are simple if and only if $\tdfB$ takes values
in $\{0,1\}$; in this case, $\sum\tdfB$ is merely the
number of triple singular points of~$\B$.
\endRemark

\section{The Zariski--van Kampen method\label{S.slopes}}

In Subsections~\ref{s.proper}--\ref{s.monodromy.infty}, we
briefly remind the classical Zariski--van Kampen approach~\cite{vanKampen}
to the computation of the fundamental group of an algebraic curve
and the construction of~\cite{degt.e6}, which makes the braid
monodromy of a genuine trigonal curve almost canonically defined.
In Subsection~\ref{s.slope}, we introduce the concept of
\emph{slope} which lets one treat a generalized trigonal curve in
terms of its trigonal model and, in particular, keep the braid
monodromy $\BG3$-valued and easily computable. In
Subsections~\ref{s.local} and~\ref{s.braid.Sk}, we compute the
local slopes and cite the results of~\cite{degt.kplets} related to
the global braid monodromy of a trigonal curve in terms of its
skeleton. Finally, in Subsection~\ref{s.quintics}, we consider a
simple example, computing the groups of irreducible
quintics.

\subsection{Proper sections and braid monodromy}\label{s.proper}
Fix a Hirzebruch surface~$\Sigma_k$, $k\ge1$,
and a genuine trigonal curve $\B\subset\Sigma_k$.
The term `section' below stands for a continuous section of (an
appropriate restriction of) the fibration $p\:\Sigma_k\to\Cp1$.

\definition\label{def.proper}
Let $\disk\subset\Cp1$ be a closed (topological) disk. A partial
section $s\:\disk\to\Sigma_k$ of~$p$ is called \emph{proper} if its
image is disjoint from both~$E$ and $\Conv\B$.
\enddefinition

The following statement is found in~\cite{degt.e6};
it is an immediate consequence of the fact
that the restriction
$p\:p\1(\disk)\sminus(E\cup\Conv\B)\to\disk$ is a locally trivial
fibration with connected fibers and contractible base.

\lemma\label{proper.homotopic}
Any disk $\disk\subset\Cp1$ admits a proper section
$s\:\disk\to\Sigma_k$.
Any two proper sections over~$\disk$ are homotopic
in the class of proper sections\rom;
furthermore, any homotopy over a fixed
point $b\in\disk$ extends to a homotopy over~$\disk$.
\qed
\endlemma

Fix a disk $\disk\subset\Cp1$ and
let $b_1,\ldots,b_r\in\disk$ be all singular and, possibly, some
nonsingular fibers of~$\B$ that belong to~$\disk$.
Denote~$F_i=p\1(b_i)$.
We assume that
all these fibers are in the interior of~$\disk$.
Denote
$\disk^\sharp=\disk\sminus\{b_1,\ldots,b_l\}$ and
fix a point
$b\in\disk^\sharp$. The restriction
$p^\sharp\:p^{-1}(\disk^\sharp)\sminus(\B\cup E)\to\disk^\sharp$
is a locally trivial fibration with a typical fiber
$F_b^\circ\sminus\B$, and any proper section $s\:\disk\to\Sigma_k$
restricts to
a section of~$p^\sharp$. Hence, given a proper
section~$s$, one can define the group
$\pi_F:=\pi_1(F^\circ_b\sminus\B,s(b))$ and the \emph{braid
monodromy} $\gm\:\pi_1(\disk^\sharp,b)\to\Aut\pi_F$.
More generally, given a path $\Gg\:[0,1]\to\disk^\sharp$ with
$\Gg(0)=b$, one can define the \emph{translation homomorphism}
$\gm_\Gg\:\pi_F\to\pi_1(F^\circ_{\Gg(1)}\sminus\B,s(b))$.

Denote by $\Gr_b\in\pi_F$ the
`counterclockwise' generator of the abelian subgroup
$\Z\cong\pi_1(F^\circ_b\sminus\Conv\B)\subset\pi_F$. (In other
words, $\Gr_b$ is the class of a large circle in~$F_b^\circ$
encompassing $\Conv\B\cap F_b^\circ$.)
Since the fibration
$p^{-1}(\disk)\sminus(\Conv\B\cup E)\to\disk$ is trivial,
hence $1$-simple, $\Gr_b$ is invariant
under the braid monodromy and is
preserved by the translation homomorphisms.
Thus, there is a canonical identification of
the elements $\Gr_{b'}$, $\Gr_{b''}$ in the fibers over any two
points $b',b''\in\disk^\sharp$; for this reason, we will omit the
subscript~$b$ in the sequel.

In this paper, we reserve the terms
`braid monodromy' and `translation homomorphism'
for the homomorphisms~$\gm$ constructed above using a
\emph{proper} section~$s$.
Under this convention,
next
lemma follows from
Lemma~\ref{proper.homotopic} and the obvious fact that the braid
monodromy is homotopy invariant.

\lemma\label{proper.braid}
The braid monodromy $\gm\:\pi_1(\disk^\sharp,b)\to\Aut\pi_F$ is well
defined and independent of the choice of a proper section
over~$\disk$ passing through $s(b)$. Given a path $\Gg$
in~$\disk^\sharp$, the translation homomorphism~$\gm_\Gg$ is
independent of the choice of a proper section passing through
$s(\Gg(0))$ and $s(\Gg(1))$ up to conjugation by~$\Gr$.
\qed
\endlemma


\subsection{The Zariski--van Kampen theorem}\label{s.vK}
Pick a basis $\{\Ga_1,\Ga_2,\Ga_3\}$ for $\pi_F$ and a basis
$\{\Gg_1,\ldots,\Gg_r\}$ for $\pi_1(\disk^\sharp,b)$.
Both $F_b^\circ\sminus\B$ and $\disk^\sharp$ are oriented
punctured planes, and we usually assume that the bases are
standard: each basis element is represented by the loop formed by
the counterclockwise boundary of a small disk centered at a
puncture and a simple arc connecting this disk to the base point;
all disks and arcs are disjoint except at the common base point.
With a certain abuse of the language, we will refer to~$\Gg_i$
(respectively,~$\Ga_j$) as the generator about the $i$-th singular
fiber (respectively, about the $j$-th branch) of~$\B$.
We also assume that the basis elements are numbered so that
$\Ga_1\Ga_2\Ga_3=\Gr$ and $\Gg_1\ldots\Gg_r$ is freely homotopic
to the boundary $\partial\disk$. Under this convention on
the basis $\{\Ga_1,\Ga_2,\Ga_3\}$, the braid monodromy does indeed
take values in the
braid group $\BG3\subset\Aut\pi_F$.

Using a proper section~$s$, we
can identify each generator~$\Gg_i$ with a certain element of the
group
$\pi_1(p^{-1}(\disk^\sharp)\sminus(\B\cup E),s(b))$;
this element does not depend on the choice of a section.
The following
presentation of the latter group
is the essence of Zariski--van Kampen's method for
computing the fundamental group of a plane algebraic curve,
see~\cite{vanKampen} for the proof and further details.

\theorem\label{vK.sharp}
In the notation above, one has
$$
\multline
\pi_1(p^{-1}(\disk^\sharp)\sminus(\B\cup E),s(b))=\\
 \bigl<\Ga_1,\Ga_2,\Ga_3,\Gg_1,\ldots,\Gg_r\bigm|
 \Gg_i\1\Ga_j\Gg_i=\gm_i(\Ga_j),\ i=1,\ldots,r,\ j=1,2,3\bigl>,
\endmultline
$$
where $\gm_i=\gm(\Gg_i)$, $i=1,\ldots,r$.
\qed
\endtheorem

\subsection{The monodromy at infinity and
relation at infinity}\label{s.monodromy.infty}
Let $\disk\subset\Cp1$ be a disk as above.
Connecting $\partial\disk$ with the base
point~$b$ by a path in~$\disk^\sharp$ and traversing it in the
counterclockwise direction (with respect to the canonical complex
orientation of~$\disk$), one obtains a certain element
$[\partial\disk]\in\pi_1(\disk^\sharp,b)$ (which depends
on the choice of the path above). The following two statements are
proved in~\cite{degt.e6}.

\lemma\label{monodromy.infty}
Assume that the interior of~$\disk$
contains all singular fibers of~$\B$. Then, for any
$\Ga\in\pi_F$, one has
$\gm([\partial\disk])(\Ga)=\Gr^{k}\Ga\Gr^{-k}$. In
particular,
the image $\gm([\partial\disk])\in\Aut\pi_F$
does not depend on the choices in
the definition of $[\partial\disk]$\rom; it is called the
\emph{monodromy at infinity}.
\qed
\endlemma

\lemma\label{relation.infty}
Assume that the interior of~$\disk$
contains all singular fibers of~$\B$. Then a presentation for
the group
$$\tsize
\pi_1(\Sigma_k\sminus(\B\cup E\cup\bigcup_{i=1}^rF_i),s(b))
$$
is obtained from that
given by Lemma~\ref{vK.sharp} by adding
the
so called \emph{relation at infinity}
$\Gg_1\ldots\Gg_r\Gr^k=1$.
\qed
\endlemma

It remains to remind that patching back in a singular fiber
$F_i$ results in an extra relation $\Gg_i=1$. Hence, for
a \emph{genuine} trigonal curve~$\B$, one has
$$
\pi_1(\Sigma_k\sminus(\B\cup E))=
 \bigl<\Ga_1,\Ga_2,\Ga_3\bigm|\gm_i=\id,\ i=1,\ldots,r,\
 \Gr^k=1\bigl>,
\eqtag\label{eq.vK}
$$
where each \emph{braid relation} $\gm_i=\id$, $i=1,\ldots,r$,
should be understood
as a triple of relations $\gm_i(\Ga_j)=\Ga_j$, $j=1,2,3$.

\subsection{Slopes}\label{s.slope}
Now, let $\BB\subset\Sigma_{\kk}$ be a
\emph{generalized} trigonal curve,
and let $\B\subset\Sigma_k$ be
its trigonal model. Denote by $F_1,\ldots,F_r$ the singular fibers
of~$\BB$ and let $b_i\in\Cp1$ be the projection of~$F_i$,
$i=1,\ldots,r$. The birational transformation between $\BB$
and~$\B$ establishes a diffeomorphism
$$\tsize
\Sigma_{\kk}\sminus(\BB\cup E\cup\bigcup_{i=1}^rF_i)\cong
\Sigma_k\sminus(\B\cup E\cup\bigcup_{i=1}^rF_i);
$$
hence, it establishes an isomorphism of the fundamental groups.
Let~$\bbranch_i$ be a small analytic disk in
$\Sigma_{\kk}\sminus E$
transversal to~$F_i$ and disjoint from~$\BB$ and from the other
singular fibers of~$\BB$, and let $\branch_i$ be the transform
of~$\bbranch_i$ in~$\Sigma_k$. We will call~$\branch_i$ a
\emph{geometric slope} of~$\BB$ at~$F_i$.
According to van Kampen's theorem~\cite{vanKampen}, patching back
in the fiber~$F_i$ results in an extra relation
$[\partial\bbranch_i]=1$ or, equivalently,
$[\partial\branch_i]=1$.

Fix a proper (with respect to the genuine trigonal curve~$\B$)
section~$s$ over a
disk $\disk\subset\Cp1$ containing the projection $p(\branch_i)$.
Pick a base point $b_i'\in p(\partial\branch_i)$ and denote
$F_i'=p\1(b_i')$
and
$\Gg_i'=[p(\partial\branch_i)]$. As above, we can
regard~$\Gg_i'$ both
as an element of $\pi_1(\disk^\sharp,b_i')$ and, \via~$s$, as an
element of
$\pi_1(p\1(\disk^\sharp)\sminus(\B\cup E),s(b_i'))$.
Furthermore, we can assume that the basis element
$\Gg_i\subset\pi_1(\disk^\sharp,b)$ introduced in
Subsection~\ref{s.vK}
has the form $\Gg_i=\zeta_i\cdot\Gg_i'\cdot\zeta_i\1$,
where $\zeta_i$ is a simple arc in~$\disk^\sharp$ connecting~$b$
to~$b_i'$.

Dragging the nonsingular fiber~$F_i'$ along $\Gg_i'$ and keeping
two points in the image of~$s$ and in~$\branch_i$, one can define
the \emph{relative braid monodromy}
$$
\gm_i^{\rel}\in\Aut\pi_1((F_i')^\circ\sminus\B,F_i'\cap\branch_i,s(b_i')).
$$

\definition\label{def.slope}
The \emph{local slope} of a generalized trigonal curve~$\BB$ at
its singular fiber~$F_i$ is the element
$\slope_i':=\gm_i^{\rel}(\xi)\cdot\xi\1\in\pi_1((F_i')^\circ\sminus\B,s(b_i'))$,
where $\xi_i$
is any path in $(F_i')^\circ\sminus\B$ connecting $s(b_i')$ and
$F_i'\cap\branch_i$. The \emph{\rom(global\rom) slope} of~$\BB$
at~$F_i$ (defined by a standard basis element~$\Gg_i$ or, more
precisely, by a path~$\zeta_i$ connecting the base point~$b$
to a point~$b_i'$
`close' to~$b_i$)
is the image
$\slope_i:=\gm_{\zeta_i}\1(\slope_i')\in\pi_F$.
\enddefinition

The following two statements are immediate consequences of the
definition.

\lemma\label{slope.defined}
The slope~$\slope_i$ is defined by the curve~$\BB$ and
generator~$\Gg_i$ up to conjugation by~$\Gr$
\rom(due to the indeterminacy of the translation homomorphism,
see Lemma~\ref{proper.braid}\rom)
and the transformation
$\slope_i\mapsto\gm_i(\Gb)\slope_i\Gb\1$, $\Gb\in\pi_F$
\rom(due to the choice of path~$\xi_i$ in the definition\rom).
\qed
\endlemma

\lemma\label{dG=slope}
In the fundamental group
$\pi_1(p\1p(\partial\branch_i)\sminus(\B\cup E),s(b_i'))$, the
conjugacy class containing $[\partial\branch_i]$ consists of all
elements of the form $\Gg_i'\slope_i'$, where $\slope_i'$ is a
local slope of~$\BB$ at~$F_i$.
\qed
\endlemma

Note that, in view of Lemma~\ref{slope.defined} and the relation
$(\Gg_i')\1\Gb\Gg_i'=\gm(\Gg_i')(\Gb)$, \cf\ Lemma~\ref{vK.sharp},
the elements $\Gg_i'\slope_i'$ do indeed form a conjugacy class.

As a consequence, in terms of the basis
$\{\Ga_1,\Ga_2,\Ga_3,\Gg_1,\ldots,\Gg_r\}$,
the relation $[\partial\branch_i]=1$ resulting from patching
the singular fiber~$F_i$
\emph{in the original surface~$\Sigma_{\kk}$}
becomes $\Gg_i=\slope_i\1$.
Eliminating~$\Gg_i$,
the
relations $\Gg_i\1\Ga_j\Gg_i=\gm_i(\Ga_j)$, $j=1,2,3$,
\cf\ Lemma~\ref{vK.sharp},
turn into the \emph{braid relations}
$$
\slope_i\Ga_j\slope_i\1=\gm_i(\Ga_j),\ j=1,2,3,
\quad\text{or}\quad\tgm_i=\id,
\eqtag\label{rel.braid}
$$
where $\tgm_i\:\Ga\mapsto \slope_i\1\gm_i(\Ga)\slope_i$ is the
\emph{twisted braid monodromy}.

Clearly, if $F_i$ is a proper fiber, then
$\branch_i=\bbranch_i$ and the path~$\xi_i$ in the definition
can be chosen so that
$\slope_i=1$. In this case $\tgm_i=\gm_i$ is still a braid.

\Remark\label{rem.slope}
In
view of
Lemma~\ref{slope.defined} and the fact that $\Gr$ is
invariant under $\gm_i$, for each fixed $i=1,\ldots,r$ the normal
subgroup of~$\pi_F$ defined by the relations
$\tgm_i=\id$ does not depend on the choice of a particular
slope~$\slope_i$, and the projection of~$\slope_i$ to the quotient
group $\pi_F/\tgm_i=\id$ is a well defined element
of this group
(depending on the curve~$\BB$ and basis element~$\Gg_i$
only). In particular, each slope commutes with~$\Gr$ (in the
corresponding quotient), making irrelevant the ambiguity in the
definition of the translation homomorphisms, see
Lemma~\ref{proper.braid}.
\endRemark

If all singular fibers are patched, hence all generators~$\Gg_i$
are eliminated,
the relation at infinity takes the form
$$
\Gr^k=\slope_r\ldots \slope_1.
\eqtag\label{rel.infinity}
$$
Finally, one obtains the following statement, \cf~\eqref{eq.vK},
expressing the fundamental
group $\pi_1(\Sigma_{\kk}\sminus(\BB\cup E))$ in
terms of the slopes and braid monodromy of the \emph{genuine}
trigonal curve~$\B$.

\corollary\label{vK.slopes}
For a generalized trigonal curve $\BB\subset\Sigma_{\kk}$ one has
$$
\pi_1(\Sigma_{\kk}\sminus(\BB\cup E))=
 \bigl<\Ga_1,\Ga_2,\Ga_3\bigm|\tgm_i=\id,\ i=1,\ldots,r,\
 \Gr^k=\slope_r\ldots \slope_1\bigl>,
$$
where each \emph{braid relation} $\tgm_i=\id$, $i=1,\ldots,r$,
should be understood
as a triple of relations $\tgm_i(\Ga_j)=\Ga_j$, $j=1,2,3$.
\qed
\endcorollary

The following statement simplifies the computation of the groups.

\proposition\label{prop.-1}
In the presentation given by Corollary~\ref{vK.slopes}, one can
omit \rom(any\rom) one of the braid relations $\tgm_i=\id$.
\endproposition

\proof
First, show that the first relation $\tgm_1=\id$ can be omitted.
Each braid relation $\tgm_i=\id$, $i=1,\ldots,r$,
can be rewritten as
$\slope_i\Ga=\gm_i(\Ga)\slope_i$, $\Ga\in\pi_F$. Hence, using all
but the first braid relations, one can rewrite
the relation at infinity~\eqref{rel.infinity}
in the form $\Gr^k=\ts_1\ldots\ts_r$,
where $\ts_i=\gm_r\circ\ldots\circ\gm_{i+1}(\slope_i)$ for
$i=1,\ldots,r-1$ and $\ts_r=\slope_r$. On the other hand, since
$\gm_r\circ\ldots\circ\gm_1$ is the conjugation by~$\Gr^{-k}$, see
Lemma~\ref{monodromy.infty}, the product
$\tgm_r\circ\ldots\circ\tgm_1$ is the conjugation by
$\Gr^{-k}\ts_1\ldots\ts_r=1$. It is the identity, and the relation
$\tgm_1=\id$ follows from
$\tgm_2=\ldots=\tgm_r=\id$.

Now, assume that the relation to be omitted is $\tgm_d=\id$ for
some $d=2,\ldots,r$. Since $\Gr$ is invariant under~$\gm_1$, the
relation $\tgm_1=\id$ implies that $\slope_1$ commutes with~$\Gr$.
Then, due to~\eqref{rel.infinity}, it also commutes with
$\slope_r\ldots\slope_2$ and \eqref{rel.infinity} is equivalent to
$\Gr^k=\slope_1\slope_r\ldots\slope_2$. Proceeding by induction
and using
only $\tgm_1=\ldots=\tgm_{d-1}=\id$,
one can rewrite~\eqref{rel.infinity} in the form
$\Gr^k=\slope_{d-1}\ldots\slope_1\slope_r\ldots\slope_d$, the
right hand side being a cyclic permutation of
$\slope_r\ldots\slope_1$. On the other hand, the cyclic permutation
$\{\Gg_d,\ldots,\Gg_r,\Gg_1,\ldots,\Gg_{d-1}\}$ is another
standard basis for $\pi_1(\disk^\sharp,b)$, with the same (but
rearranged)
slopes~$\slope_i$
and braid monodromies~$\gm_i$,
and with respect to this new basis
the braid relation to be omitted
is the first one. Hence the statement follows from the first part of
the proof.
\endproof

\Remark\label{rem.cyclic.order}
A by-product of the previous proof is the fact that, modulo (all
but one) braid relations the slopes
$\slope_r,\ldots,\slope_1$ cyclically commute, \ie,
one has
$$
\slope_r\ldots \slope_2\slope_1=
\slope_{d-1}\ldots\slope_1\slope_r\ldots\slope_d
$$
for each $d=2,\ldots,r$.
In particular, if there are only two nontrivial slopes (which is
always the case in Section~\ref{S.proof.group} below), their order
in the relation at infinity~\eqref{rel.infinity} is irrelevant.
\endRemark

\subsection{Local slopes and braid relations}\label{s.local}
Given two
elements~$\Ga$, $\Gb$ of a group and a nonnegative integer~$m$,
introduce the notation
$$
\{\Ga,\Gb\}_m=\cases
(\Ga\Gb)^k(\Gb\Ga)^{-k},&\text{if $m=2k$ is even},\\
\bigl((\Ga\Gb)^k\Ga\bigr)\bigl((\Gb\Ga)^k\Gb\bigr)\1,&
 \text{if $m=2k+1$ is odd}.
\endcases
$$
The relation $\{\Ga,\Gb\}_m=1$ is equivalent to $\Gs^m=\id$, where
$\Gs$ is the Artin generator of the braid group~$\BG2$ acting on the
free group $\<\Ga,\Gb\>$. Hence,
$$
\{\Ga,\Gb\}_m=\{\Ga,\Gb\}_n=1\quad
\text{is equivalent to}\quad
\{\Ga,\Gb\}_{\gcd(m,n)}=1.
\eqtag\label{eq.equiv}
$$
For the small values of~$m$, the relation $\{\Ga,\Gb\}_m=1$ takes
the following form:
\Dashes
\dash
$m=0$: tautology;
\dash
$m=1$: the identification $\Ga=\Gb$;
\dash
$m=2$: the commutativity relation $[\Ga,\Gb]=1$;
\dash
$m=3$: the braid relation $\Ga\Gb\Ga=\Gb\Ga\Gb$.
\endDashes

Let $F_i$ be a type $\tA{_p}$ (type $\tA{_0^*}$ if $p=0$) singular
fiber of a trigonal curve~$B$, and let $b_i=p(F_i)\subset\disk$ be
its projection. Pick a simple arc $\zeta_i\:[0,1]\to\disk$
connecting the base point~$b$ to~$b_i$ and such that
$\zeta_i([0,1))\subset\disk^\sharp$. We say that two
consecutive generators~$\Ga_j$, $\Ga_{j+1}$ of~$\pi_F$,
$j=1$ or~$2$,
\emph{collide at~$F_i$} (along~$\zeta_i$) if there is a Milnor
ball~$\MB$ about the point of non-transversal intersection of~$\B$
and~$F_i$ such that, for each sufficiently small $\epsilon>0$, the
images of~$\Ga_j$ and~$\Ga_{j+1}$ under the translation
along the restriction $\zeta_i|_{[0,1-\epsilon]}$ are
represented by a pair of loops that differ only inside~$\MB$, whereas the
image of the third generator is represented by a loop totally
outside~$\MB$.

In~\ref{E.proper}--\ref{E.typeA} below, we pick a standard basis
$\{\Ga_1,\Ga_2,\Ga_3\}$ for~$\pi_F$ and
assume that two consecutive
elements of this basis  collide at~$F_i$ along a certain
path~$\zeta_i$; then we use
this path~$\zeta_i$ to construct a generator~$\Gg_i$
about~$F_i$
as in
Subsection~\ref{s.vK}. In other words, it is $\zeta_i$ that is
used to define the global braid monodromy~$\gm_i$ and the
slope~$\slope_i$. All computations are straightforward, using
local normal forms of the singularities involved; we merely state
the results.

We denote by~$\Gs_1$, $\Gs_2$ the Artin generators of the braid
group~$\BG3$ acting on the free group
$\pi_F=\<\Ga_1,\Ga_2,\Ga_3\>$, so that
$\Gs_i\:\Ga_i\mapsto\Ga_i\Ga_{i+1}\Ga_i\1$,
$\Ga_{i+1}\mapsto\Ga_i$, $i=1,2$.

Note that the slope~$\slope_i$ is only useful if the fiber~$F_i$
is to be patched back in (as otherwise the presentation for the
group $\pi_1(\Sigma_{\kk}\sminus(\BB\cup E\cup F_i\cup\ldots))$
would contain the original generator~$\Gg_i$ rather
than~$\slope_i$). For this reason, after a small equisingular
deformation of $\BB+E$, one can assume that $\BB$ is maximally
transversal to~$F_i$ at infinity. We do make this assumption below.

\subsubsection{A proper fiber}\label{E.proper}
Assume that $F_i$ is a proper type~$\tA{_p}$
(type~$\tA{_0^*}$ if $p=0$) fiber of $\BB=\B$. Then one has:
\roster
\item\local{proper.12}
if $\Ga_1$ and $\Ga_2$
collide at~$F_i$,
then
$\gm_i=\Gs_1^{p+1}$ and $\slope_i=1$, so that
the braid relations are
$\{\Ga_1,\Ga_2\}_{p+1}=1$;
\item\local{proper.23}
if $\Ga_2$ and $\Ga_3$
collide at~$F_i$,
then
$\gm_i=\Gs_2^{p+1}$ and $\slope_i=1$, so that
the braid relations are
$\{\Ga_2,\Ga_3\}_{p+1}=1$.
\endroster

\subsubsection{A nonsingular branch at infinity}\label{E.smooth}
Assume that $F_i$ is a type $\tA{_p}$ singular fiber
(type~$\tA{_0^*}$ if $p=0$ or no singularity if $p=-1$) and
a single smooth branch of~$\BB$
intersects~$E$ at $F_i\cap E$ with multiplicity $q\ge1$. Then
$F_i$ is a type $\tA{_{p+2q}}$ singular fiber of~$\B$
and
one has:
\roster
\item\local{smooth.12}
if $\Ga_1$ and $\Ga_2$
collide at~$F_i$,
then
$\gm_i=\Gs_1^{p+2q+1}$ and $\slope_i=(\Ga_1\Ga_2)^q$, so that
the braid relations are
$\{\Ga_1,\Ga_2\}_{p+1}=[(\Ga_1\Ga_2)^q,\Ga_3]=1$;
\item\local{smooth.23}
if $\Ga_2$ and $\Ga_3$
collide at~$F_i$,
then
$\gm_i=\Gs_2^{p+2q+1}$ and $\slope_i=(\Ga_2\Ga_3)^q$, so that
the braid relations are
$\{\Ga_2,\Ga_3\}_{p+1}=[\Ga_1,(\Ga_2\Ga_3)^q]=1$.
\endroster

\subsubsection{A double point at infinity}\label{E.typeA}
Assume that $\BB$ has a type $\bA_p$ singular point at $F_i\cap E$
and intersects~$E$ at this point with multiplicity $2q$,
$1\le q\le(p+1)/2$. Then $F_i$ is a type $\tA{_{p-2q}}$ singular
fiber of~$\B$ (with the same convention as
in~\ref{E.smooth} for
the values $p-2q=0$ or $-1$)
and
one has:
\roster
\item\local{typeA.12}
if $\Ga_1$ and $\Ga_2$
collide at~$F_i$,
then
$\gm_i=\Gs_1^{p-2q+1}$ and $\slope_i=\Ga_3^q$, so that
the braid relations are
$\Gs_1^{p-2q+1}(\Ga_j)=\Ga_3^q\Ga_j\Ga_3^{-q}$, $j=1,2$;
\item\local{typeA.23}
if $\Ga_2$ and $\Ga_3$
collide at~$F_i$,
then
$\gm_i=\Gs_2^{p-2q+1}$ and $\slope_i=\Ga_1^q$, so that
the braid relations are
$\Gs_2^{p-2q+1}(\Ga_j)=\Ga_1^q\Ga_j\Ga_1^{-q}$, $j=2,3$.
\endroster
(If $p-2q=-1$, there is no collision and $\gm_i=\id$. In this
case, the slope is $\slope_i=\Ga_j^q$, where $\Ga_j$ is the
generator about the proper branch of the original curve~$\BB$.)

Now, assume that $\BB$ has a type $\bA_{2p-1}$ singular point at
$F_i\cap E$ and the two branches at this point intersect~$E$ with
multiplicities~$p$ and~$p+q$ for some $q\ge1$. Then $F_i$ is a
type $\tA{_{2q-1}}$ singular fiber of~$\B$, and one of the two
branches of~$\B$ at its type $\bA_{2q-1}$ singular point in~$F_i$
is distinguished: it is the transform of the proper
branch
of~$\BB$. Choose generators $\Ga_1$, $\Ga_2$, $\Ga_3$ so that
$\Ga_1$ and~$\Ga_2$ collide at~$F_i$ and~$\Ga_1$ is the generator
about the distinguished branch of~$\B$. Then
\roster
\item[3]
$\gm_i=\Gs_1^{2q}$ and $\slope_i=(\Ga_1\Ga_2)^q\Ga_1^p$,
so that
the braid relations are
$[\Ga_1^p,\Ga_2]=1$ and
$[(\Ga_1\Ga_2)^q\Ga_1^p,\Ga_3]=1$.
\endroster

Finally, assume that $\BB$ has a type~$\bA_{2p}$ singular point
at $F_i\cap E$ and intersects~$E$ at this point with
multiplicity $2p+1$. Then $F_i$ is a type~$\tA{_1^*}$ singular
fiber of~$\B$ and, in an appropriate basis
$\{\Ga_1,\Ga_2,\Ga_3\}$, such that $\Ga_2$ is the generator about
the branch of~$\B$ transversal to~$F_i$, one has
\roster
\item[4]
$\gm_i=\Gs_1\Gs_2\Gs_1$ and $\slope_i=\Ga_1\Ga_2^{p+1}$,
so that
the braid relations are
$[\Ga_1,\Ga_2^{p+1}]=1$ and
$\Ga_3=\Ga_2^p\Ga_1\Ga_2^{-p}$.
\endroster

\subsubsection{A triple point at infinity}\label{E.triple}
Assume that $\BB$ has a triple point at $F_i\cap E$ and consider
the (generalized) trigonal curve~$\BB'$ obtained from~$\BB$ by one
elementary transformation centered at this point. Then the
transform~$\bbranch'_i$ of~$\bbranch_i$
is still disjoint from~$\BB'$
and thus can be used to define the slope of~$\BB$;
hence the slope of~$\BB$ at~$F_i$ equals that of~$\BB'$. As a
consequence, one has the following statement.

\corollary
For $\BB\subset\Sigma_{\kk}$ and $\B''\subset\Sigma_{\kk+1}$ as
above there is a canonical isomorphism
$\pi_1(\Sigma_{\kk}\sminus(\BB\cup E))=
\pi_1(\Sigma_{\kk+1}\sminus(\BB'\cup E))$.
\qed
\endcorollary

\subsection{Braid monodromy via skeletons}\label{s.braid.Sk}
Let $\B\subset\Sigma_k$ be a trigonal curve, and let
$\Sk\subset\Cp1$ be its skeleton. Below, we cite a few results
of~\cite{degt.kplets} concerning the braid monodromy of~$\B$ in terms
of~$\Sk$. For simplicity, we assume that all \black-vertices
of~$\Sk$ are trivalent and all its \white-vertices are bivalent
(hence omitted). An alternative description, including more
general skeletons, is found in~\cite{tripods}.

Recall that a \emph{marking} at a trivalent \black-vertex~$v$
of~$\Sk$ is a counterclockwise order $e_1$, $e_2$, $e_3$
of the three edges adjacent to~$v$, see Figure~\ref{fig.basis}(a).
We consider the indices
defined modulo~$3$, so that $e_{i+3}=e_i$. The three points of
intersection of~$\B$ and the fiber~$F_v$ over~$v$ form an
equilateral triangle. These points are in a canonical one-to-one
correspondence with the edges~$e_i$ of~$\Sk$ at~$v$. Hence, a
marking gives rise to a \emph{canonical basis}
$\{\Ga_1,\Ga_2,\Ga_3\}$ of the group
$\pi_F=\pi_1(F^\circ_v\sminus\B)$, see
Figure~\ref{fig.basis}(b); this basis is well defined up to
simultaneous conjugation by a power of $\Gr=\Ga_1\Ga_2\Ga_3$.

\midinsert
\centerline{\vbox{\halign{\hss#\hss&&\qquad\qquad\hss#\hss\cr
\cpic{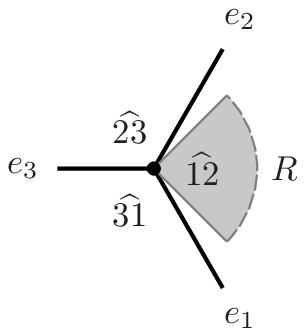}&
\cpic{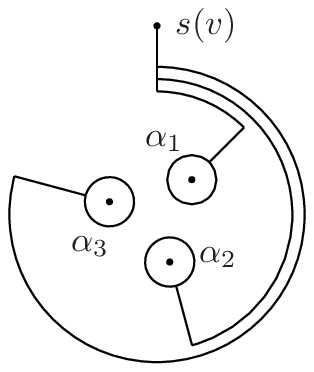}\cr
\noalign{\medskip}
(a)&(b)\cr
\crcr}}}
\figure\label{fig.basis}
A marking (a) and a canonical basis (b)
\endfigure
\endinsert

As in Subsection~\ref{s.local}, let $\Gs_1$ and~$\Gs_2$ be the
Artin generators of the braid group~$\BG3$ acting on
$\<\Ga_1,\Ga_2,\Ga_3\>$. Denote $\Gs_3=\Gs_1\1\Gs_2\Gs_1$ and
extend indices to all integers
\via\ $\Gs_{i\pm3}=\Gs_i$. Note that the map
$(\Gs_{i-1},\Gs_i)\mapsto(\Gs_i,\Gs_{i+1})$ is an
automorphism of~$\BG3$.
Recall also that the center of~$\BG3$ is the cyclic subgroup
generated by
$(\Gs_1\Gs_2)^3=(\Gs_2\Gs_3)^3=(\Gs_3\Gs_1)^3$.

\subsubsection{The translation homomorphisms}\label{Sk.tr}
Let~$u$ and~$v$ be two marked
\black-vertices of~$\Sk$
connected by a single edge~$e$; to indicate the markings, we
use the notation $e=[i,j]$, where $i$ and~$j$ are the indices
of~$e$ at~$u$ and~$v$, respectively. Choosing a pair of
canonical bases
defined by the markings, one can identify the groups
$\pi_1(F^\circ_u\sminus\B)$ and
$\pi_1(F^\circ_v\sminus\B)$ with the `standard' free group
$\<\Ga_1,\Ga_2,\Ga_3\>$ and thus regard the translation
homomorphism
$\gm_e\:\pi_1(F^\circ_u\sminus\B)\to\pi_1(F^\circ_v\sminus\B)$
as an automorphism
of $\<\Ga_1,\Ga_2,\Ga_3\>$. It is a braid. However, since both the
bases and the homomorphism~$\gm_e$ itself are only defined up to
conjugation by~$\Gr$ (unless a proper section is fixed, see
Lemma~\ref{proper.braid}), this automorphism should be regarded as
an element of the reduced braid group
$\bar\BG3=\BG3/(\Gs_1\Gs_2)^3\cong\PSL(2,\Z)$.
On the other hand,
this ambiguity does not affect the computation
of the fundamental group, \cf\ Remark~\ref{rem.slope}.

With the above convention,
the translation homomorphism $\gm_{[i,j]}\in\bar\BG3$
along an edge $e=[i,j]$
is given as follows:
$$
\gm_{[i,i+1]}=\Gs_i,\quad
\gm_{[i+1,i]}=\Gs_i\1,\quad
\gm_{[i,i]}=\Gs_i\Gs_{i-1}\Gs_i,\quad
i\in\Z.
$$
The translation homomorphism $\gm_\Gg\in\bar\BG3$
along a path~$\Gg$ composed
by edges of~$\Sk$ is the composition of the contributions of
single
edges. If $\Gg$ is a loop, the braid monodromy $\gm_\Gg$ is a well
defined element of~$\BG3$. It is uniquely recovered from its
projection to $\bar\BG3$ just described and its \emph{degree}
(\ie, the image in $\BG3/[\BG3,\BG3]=\Z$); the latter is
determined by the number and the types of the singular fibers
of~$\B$ encompassed by~$\Gg$. More precisely,
for a disk $\disk\subset\Cp1$ as in Subsection~\ref{s.proper},
the composed
homomorphism $\pi_1(\disk^\sharp)\to\BG3\to\Z$ sends a
generator~$\Gg_i$ about a singular fiber~$F_i$ to the
multiplicity~$\mult F_i$, see Subsection~\ref{s.trigonal}.

\subsubsection{The braid relations resulting from a region}\label{Sk.m}
Given a trivalent \black-vertex~$v$ of~$\Sk$, one can define three
(germs of) angles at~$v$, which are represented by the connected
components of the intersection of $\Cp1\sminus\Sk$ and a regular
neighborhood of~$v$ in~$\Cp1$. If $v$ is marked, we denote these
angles $\corner12$, $\corner23$, and $\corner31$, according to the
two edges adjacent to an angle, see Figure~\ref{fig.basis}(a).
The position of a region~$R$ adjacent to~$v$ with respect to the
marking at~$v$ can then be described by indicating the angle(s)
that belong to~$R$; for example, in Figure~\ref{fig.basis}(a)
one has $\corner12\subset R$. Note that a region may contain two
or even all three angles at~$v$,
see \eg\ the outer nonagon and the central vertex
in Figure~\ref{fig.d}(b) below.

Let $\disk\subset\Cp1$ be a closed disk as in
Subsection~\ref{s.proper}. Assume that $v\in\partial\disk$ and
that $\disk\sminus v\subset R$, \cf\ the shaded area in
Figure~\ref{fig.basis}(a). Then $\disk$ intersects exactly one of
the three angles at~$v$ (in the figure this angle is~$\corner12$).
Take $v=b$ for the base point and let $\{\Ga_1,\Ga_2,\Ga_3\}$
be a canonical basis for $\pi_F=\pi_1(F^\circ_v\sminus\B)$
defined by the marking
at~$v$. The following three statements are straightforward; for
details see~\cite{degt.kplets}.

\lemma\label{lem.collide}
If a disk~$\disk$ as above intersects angle~$\corner12$
\rom(respectively,~$\corner23$\rom), then $\Ga_1$ and~$\Ga_2$
\rom(respectively, $\Ga_2$ and~$\Ga_3$\rom) collide at any
type~$\tA{}$ singular fiber of~$\B$ in~$\disk$ along any path
contained in~$\disk$.
\qed
\endlemma

\lemma\label{abelian.monodromy}
If a disk~$\disk$ as above intersects angle~$\corner12$
\rom(respectively,~$\corner23$ or~$\corner31$\rom),
the braid monodromy $\gm\:\pi_1(\disk^\sharp,v)\to\Aut\pi_F$ takes
values in the abelian subgroup of $\BG3\subset\Aut\pi_F$ generated
by
the central element $(\Gs_1\Gs_2)^3$ and~$\Gs_1$ \rom(respectively,
$\Gs_2$ or~$\Gs_3$\rom).
If all singular fibers in~$\disk$ are of type~$\tA{}$, then $\gm$
takes values in the cyclic subgroup generated by~$\Gs_1$
\rom(respectively, by~$\Gs_2$ or~$\Gs_3$\rom).
\qed
\endlemma

More precisely, in Lemma~\ref{abelian.monodromy}, the value
of~$\gm$ on a type~$\tA{_{p-1}}$ fiber (type~$\tA{_0^*}$ if $p=1$)
is $\Gs_i^p$ (assuming that $\disk$ intersects the angle spanned
by~$e_i$ and~$e_{i+1}$), and its value on a type~$\tD{_{q+4}}$
fiber is $(\Gs_1\Gs_2)^3\Gs_i^q$. The value at a
non-simple singular fiber of type $\tJ{_{r,p}}$
is $(\Gs_1\Gs_2)^{3r}\Gs_i^p$.

\corollary\label{region.relations}
Assume that a region~$R$ of~$\Sk$
adjacent to a marked vertex~$v$
contains, among others, singular fibers of types
$\tA{_{p_i-1}}$ \rom($\tA{_0^*}$ if
$p_i=1$\rom), $i=1,\ldots,s$. Denote $p=\gcd(p_i)$.
Then the braid
relations $\gm_i=\id$
resulting from these fibers are equivalent to a single
relation as follows\rom:
\Dashes
\dash
$\{\Ga_1,\Ga_2\}_p=1$ if $\corner12\subset R$\rom;
\dash
$\{\Ga_2,\Ga_3\}_p=1$ if $\corner23\subset R$\rom;
\dash
$\{\Ga_1,\Ga_2\Ga_3\Ga_2\1\}_p=1$ if $\corner31\subset R$.
\endDashes
In particular, if an $m$-gonal region~$R$
contains a single singular
fiber, which is of type~$\tA{}$, it results in a single braid
relation as above with $p=m$.
\qed
\endcorollary

\Remark
In Corollary~\ref{region.relations}, if $\B$ is the trigonal model of
a generalized trigonal curve $\BB\subset\Sigma_{\kk}$
and it is $\pi_1(\Sigma_{\kk}\sminus(\BB\cup E))$ that is
computed, one should assume in addition that the fibers considered
are proper for~$\BB$.
\endRemark

\subsubsection{An irreducibility criterion}
A \emph{marking} of a
skeleton~$\Sk$ is a collection of markings at all its trivalent
\black-vertices. A
marking of a generic skeleton without singular \black-vertices
is called \emph{splitting} if it satisfies the following
three conditions:
\roster
\item\local{splitting.1}
the types of all edges, \cf~\ref{Sk.tr},
are $[1,1]$, $[2,3]$, or~$[3,2]$;
\item\local{splitting.2}
an edge connecting a \black-vertex~$v$ and a singular
\white-vertex has index~$1$ at~$v$;
\item\local{splitting.3}
if a region~$R$ contains angle $\corner12$ or $\corner31$ at one
of its vertices, the multiplicities of all singular fibers
inside~$R$ are even.
\endroster
(Note that, given~\loccit{splitting.1} and~\loccit{splitting.2},
the last condition holds automatically if $R$ contains a single
singular fiber, as $R$ is necessarily a $(2m)$-gon.) The following
criterion is essentially contained in~\cite{degt.kplets}; it is
obtained by reducing the braid monodromy to the symmetric
group~$\SG3$.

\theorem\label{th.splitting}
A trigonal curve $\B\subset\Sigma_k$ with connected generic
skeleton~$\SkB$
is reducible if and only if\/ $\SkB$
has no singular \black-vertices and admits a splitting marking.
Each such marking defines a component of~$\B$
that is a section of~$\Sigma_k$.
\qed
\endtheorem

\Remark\label{rem.splitting}
A splitting marking defines a component of~$\B$ as follows: over
each \black-vertex~$v$, in a canonical basis
$\{\Ga_1,\Ga_2,\Ga_3\}$ defined by the marking, $\Ga_1$ is the
generator about the distinguished component.
\endRemark

\subsection{Example: irreducible quintics}\label{s.quintics}
As a simple example of application of the techniques developed
in this section, we recompute the non-abelian fundamental groups
of irreducible plane quintics, see~\cite{groups} and~\cite{Artal}.
A more advanced example is the
contents of Section~\ref{S.proof.group} below.

Let $\CC\subset\Cp2$ be a quintic with the set of singularities
$\bA_6\splus3\bA_2$ or $3\bA_4$. Blow up the type~$\bA_6$ point
(respectively, one of the type~$\bA_4$ points) to obtain a
generalized trigonal curve $\BB\subset\Sigma_1$ and let
$\B\subset\Sigma_2$ be the trigonal model of~$\BB$. It is a
maximal trigonal curve with the
combinatorial type of
singular fibers $4\tA{_2}$ or
$2\tA{_4}\splus2\tA{_0^*}$; its skeleton~$\Sk$ is shown in
Figures~\ref{fig.quintics}(a) and~(b), respectively.

\midinsert
\centerline{\vbox{\halign{\hss#\hss&&\qquad\qquad\hss#\hss\cr
\cpic{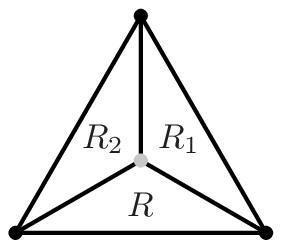}&
\cpic{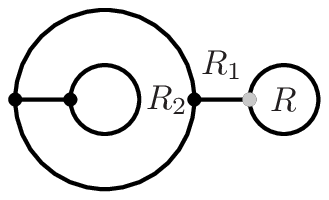}\cr
\noalign{\medskip}
(a)&(b)\cr
\crcr}}}
\figure\label{fig.quintics}
Skeletons of plane quintics
\endfigure
\endinsert

Let~$R$ be the region of~$\Sk$ containing the only improper fiber
of~$\BB$. We choose for the reference point~$b$ the vertex shown
in the figures in grey and take for $\{\Ga_1,\Ga_2,\Ga_3\}$ a
canonical basis over~$b$ defined by the marking such that
$\corner12\subset R$.
In both cases, the only nontrivial slope $\slope=\Ga_3$
is given
by~\iref{E.typeA}{typeA.12}, with $(p,q)=(4,1)$ or $(2,1)$,
respectively. According to Proposition~\ref{prop.-1}, the
fundamental group
$$
\pi_1:=\pi_1(\Cp2\sminus\CC)=\pi_1(\Sigma_1\sminus(\BB\cup E))
$$
is defined by the relation at infinity and the braid relations
resulting from three (out of four) regions~$R$, $R_1$, $R_2$ shown
in the figures. Using Subsection~\ref{s.braid.Sk}, one obtains the
following relations: for the set of singularities
$\bA_6\splus3\bA_2$, see Figure~\ref{fig.quintics}(a):
$$
\gather
\Gr^2=\Ga_3,\quad
 \{\Ga_2,\Ga_3\}_3=\{\Ga_1,\Ga_2\Ga_3\Ga_2\1\}_3=1,\\
(\Ga_1\Ga_2\Ga_1)\Ga_2(\Ga_1\Ga_2\Ga_1)\1=\Ga_3\Ga_1\Ga_3\1,\quad
 (\Ga_1\Ga_2)\Ga_1(\Ga_1\Ga_2)\1=\Ga_3\Ga_2\Ga_3\1,
\endgather
$$
and for the set of singularities $3\bA_4$, see
Figure~\ref{fig.quintics}(b):
$$
\Gr^2=\Ga_3,\quad
 \{\Ga_2,\Ga_3\}_5=\{\Ga_2,\Gr\1\Ga_1\Gr\}_5,\quad
\Ga_1\Ga_2\Ga_1\1=\Ga_3\Ga_1\Ga_3\1,\quad
 \Ga_1=\Ga_3\Ga_2\Ga_3\1.
$$
In the former case, the group is known to be infinite,
see~\cite{groups}, as it factors to infinite
Coxeter's group $(2,3,7)$,
see~\cite{Coxeter}.
In the latter case, using \GAP~\cite{GAP}, one can see that
$\pi_1$ is a soluble group of order~$320$; one has
$$
\pi_1/\pi_1'=\CG5,\quad
\pi_1'/\pi_1''=(\CG2)^4,\quad
\pi_1''=(\CG2)^2,
$$
where $'$ stands for the commutant.

Certainly, this approach applies as well to other quintics with a
double singular point. In particular, one can easily show that the
groups of all other irreducible quintics are abelian.

\section{Proof of Theorem~\ref{th.d}\label{S.proof.d}}

In Subsection~\ref{s.models}, we replace a
sextic as in the theorem with its \emph{trigonal model},
which is a maximal trigonal curve in an
appropriate Hirzebruch surface. Then, in Subsections~\ref{s.D7}
and~\ref{s.D5}, we enumerate the possible
skeletons of trigonal models; in view of
Theorem~\ref{th.skeleton}, this enumeration suffices to prove
Theorem~\ref{th.d}.

\subsection{The trigonal models}\label{s.models}
Recall that, due to~\cite{dessin.e7},
an irreducible maximizing simple sextic cannot have a singular
point of type~$\bD_{2k}$, $k\ge2$ or more than one singular point
from the list
$\bA_{2k+1}$, $k\ge0$, $\bD_{2k+1}$, $k\ge2$, or~$\bE_7$. Thus, a
sextic as in Theorem~\ref{th.d} has a unique type~$\bD$ point,
which is either~$\bD_5$ or $\bD_p$ with odd $p\ge7$. We will
consider the two cases separately.

Both Theorem~\ref{th.d} and Theorem~\ref{th.group} are proved by a
reduction of sextics to trigonal curves. A key r\^ole is played by
the following two propositions.

\proposition\label{1-1.d}
There is a natural
bijection~$\phi$, invariant under equisingular deformations,
between Zariski open and dense in each equisingular stratum
subsets of the following two sets\rom:
\roster
\item\local{d.sextic}
plane sextics~$\CC$ with a distinguished
type~$\bD_p$, $p\ge7$,
singular point~$P$ and without linear components
through~$P$, and
\item\local{d.trigonal}
trigonal curves $\B\subset\Sigma_3$ with a distinguished
type~$\tA{_1}$ singular fiber $\Fi$ and a distinguished
type~$\tA{_{p-7}}$ \rom($\tA{_0^*}$ if $p=7$\rom)
singular fiber $\Fii\ne \Fi$.
\endroster
A sextic~$\CC$ is irreducible if and only if so is $\B=\phi(\CC)$, and
$\CC$ is maximizing if and only if $\B$ is
maximal and stable.
\endproposition

\proposition\label{1-1.d5}
There is a natural
bijection~$\phi$, invariant under equisingular deformations,
between Zariski open and dense in each equisingular stratum
subsets of the following two sets\rom:
\roster
\item\local{d.sextic}
plane sextics~$\CC$ with a distinguished
type~$\bD_5$
singular point~$P$ and without linear components
through~$P$, and
\item\local{d.trigonal}
trigonal curves $\B\subset\Sigma_4$ with a distinguished
type~$\tA{_1}$ singular fiber $\Fi$ and a distinguished
type~$\tA{_3}$ singular fiber~$\Fii$.
\endroster
A sextic~$\CC$ is irreducible if and only if so is $\B=\phi(\CC)$, and
$\CC$ is maximizing if and only if $\B$ is
maximal and stable.
\endproposition

The trigonal curve $\B=\phi(\CC)$ corresponding to a sextic~$\CC$
\via\ Propositions~\ref{1-1.d} and~\ref{1-1.d5} is called the
\emph{trigonal model} of~$\CC$.

\proof[Proof of Propositions~\ref{1-1.d} and~\ref{1-1.d5}]
Let $\CC\subset\Cp2$ and~$P$ be a pair as in the statement. Blow~$P$
up and denote by $\BB\subset\Sigma_1=\Cp2(P)$
the proper transform of~$\CC$; it is a generalized trigonal curve
with two points at infinity. We let $\B=\phi(\CC)$ to be the
trigonal model of~$\BB$. The inverse transformation consists in
the passage from~$\B$ back to $\BB\subset\Sigma_1$ and blowing
down  the exceptional section of~$\Sigma_1$.

The distinguished fibers~$\Fi$ and~$\Fii$ of~$\B$ correspond to the
two points of~$\BB$ at infinity ($\Fi$ corresponding to the smooth
branch of~$\CC$ at~$P$). It is straightforward that, generically,
the types of these fibers are as indicated in the statements. If
the original sextic~$\CC$ is in a special position with respect to
the pencil of lines through~$P$, these fibers may degenerate:
$\Fi$ may degenerate to~$\tA{_1^*}$
or~$\tA{_s}$, $s>1$, and $\Fii$ may degenerate to $\tA{_0^{**}}$ (in
Proposition~\ref{1-1.d} with $p=7$) or to~$\tA{_s}$, $s>3$ (in
Proposition~\ref{1-1.d5}). However, using theory of trigonal
curves (perturbations of dessins), one can easily see that any
such curve~$\B$ can be perturbed to a generic one, and this
perturbation is followed by \emph{equisingular} deformations
of~$\BB$ and~$\CC$.

Since $\CC$ and~$\B$ are birational transforms of each other, they
are either both reducible or both irreducible. The fact that
maximizing sextics correspond to stable maximal trigonal curves
follows from
Theorem~\ref{th.maximal}:
for a generic sextic~$\CC$ as in the statements, one has
$\mu(\B)=\mu(\CC)-6$ in Proposition~\ref{1-1.d} and
$\mu(\B)=\mu(\CC)-1$ in Proposition~\ref{1-1.d5}.
\endproof

Let~$\CC$ be a sextic as in Theorem~\ref{th.d}. Since $\CC$ is
irreducible, it has a unique type~$\bD$ point (see above), which
we take for the distinguished point~$P$.
Denote by
$\B\subset\Sigma_k$, $k=3$ or~$4$, the trigonal model of~$\CC$ and
let~$\Sk$ be the skeleton of~$\B$. Let, further, $\Ri$ and~$\Rii$
be the regions of~$\Sk$ containing the distinguished singular
fibers $\Fi$ and~$\Fii$, respectively. Since we assume that $\CC$
has no type~$\bE$ singular points, $\B$ has no triple points (one
has $\tdfB\equiv0$, see Remark~\ref{rem.td})
and
hence $\Sk$ has exactly $2k$ \black-vertices and has no singular
vertices. Thus, due to Theorem~\ref{th.skeleton},
the proof of Theorem~\ref{th.d} reduces to the
enumeration of $3$-regular skeletons of irreducible curves with a
prescribed number of vertices and with a pair
$(\Ri,\Rii)$ of distinguished regions. This is done is
Subsections~\ref{s.D7} and~\ref{s.D5} below.

\subsection{The case $\bD_p$, $p\ge7$}\label{s.D7}
In this case, $\Sk$ has six \black-vertices, $\Ri$ is a bigon, and
$\Rii$ is a $(p-6)$-gon. Since $p$ is not fixed, one can take
for~$\Rii$ any region of~$\Sk$ other than~$\Ri$.

\midinsert
\centerline{\vbox{\halign{\hss#\hss&&\qquad\qquad\hss#\hss\cr
\cpic{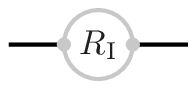}&
\cpic{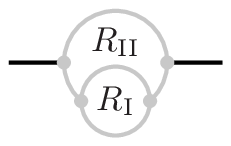}\cr
\noalign{\medskip}
(a)&(b)\cr
\crcr}}}
\figure\label{fig.insertion}
A bigonal (a) and bibigonal (b) insertions
\endfigure
\endinsert

The bigonal region~$\Ri$ looks as shown in grey in
Figure~\ref{fig.insertion}(a); we will call this region
the \emph{insertion}.
Removing~$\Ri$ from~$\Sk$ and patching the two black edges in the
figure to a single edge results in a new $3$-regular
skeleton~$\Sk'$ with four \black-vertices.
Conversely, starting from~$\Sk'$ and placing an insertion at the
middle of any of its edges produces a skeleton~$\Sk$ with six
\black-vertices and a distinguished bigonal region, which we take
for~$\Ri$. Using Theorem~\ref{th.splitting}, one can see
that $\Sk$ is the skeleton of an irreducible trigonal curve if and
only if so is~$\Sk'$. There are three such skeletons (see
\eg~\cite{symmetric}); they are listed in Figure~\ref{fig.d}.
Starting from one of these skeletons and varying the position of
the insertion (shown in grey and numbered in the figure)
and the choice of the
second distinguished region~$\Rii$, all up to symmetries of the
skeleton,
one obtains the $22$ deformation families listed in
Table~\ref{tab.d}. (Some rows of the table represent pairs of
complex conjugate
curves,
see comments below.)

\midinsert
\centerline{\vbox{\halign{\hss#\hss&&\qquad\hss#\hss\cr
\cpic{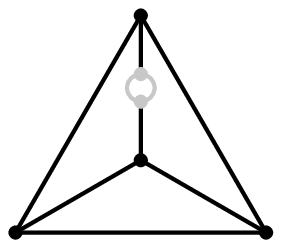}&
\cpic{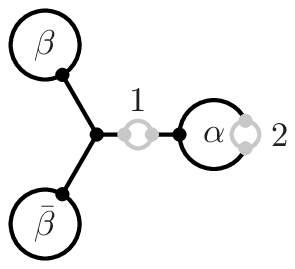}&
\cpic{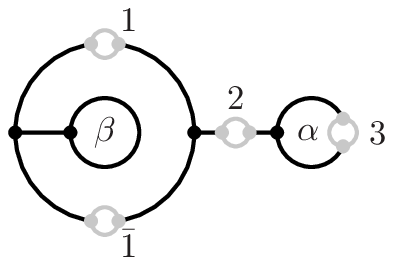}\cr
\noalign{\medskip}
(a)&(b)&(c)\cr
\crcr}}}
\figure\label{fig.d}
A type $\bD_p$ singular point, $p\ge7$
\endfigure
\endinsert

\Remark[Comments to Tables~\ref{tab.d} and~\ref{tab.d5}]
Listed in the tables are combinatorial types of singularities
and references to
the figures representing the corresponding skeletons. Equal
superscripts precede combinatorial types shared by several items
in the tables. The `Count' column lists the numbers $(n_r,n_c)$ of
real curves and pairs of complex conjugate curves, so that the
total number of curves represented by a row is
$n_r+2n_c$. The last two columns refer to the computation of the
fundamental group and indicate the parameters used in this
computation. (A parameter list is marked with a~$^*$ when the
general approach does not work quite well for a particular curve.
In this case, more details are found in the subsection referred to
in the table.)
\endRemark

\midinsert
\table\label{tab.d}
Maximal sets of singularities with a type~$\bD_p$ point, $p\ge7$
\endtable

\def\no{}
\def\*{\llap{$^*$}}
\def\same{\afterassignment\dosame\count0=}
\def\dosame{\llap{$^{\the\count0\,}$}}
\def\tref#1{\,\text{\ref{#1}}\,}
\centerline{\vbox{\offinterlineskip\halign{%
\tabstrut\ \hss#\hss\ \vrule
&\quad$#$\hss\quad\vrule
&\quad\null#\hss\quad\vrule
&&\ \hss$#$\hss\ \vrule\cr
\noalign{\hrule}
\exstrut&&&&&\cr
\#&\text{Set of singularities}&
 \hss Figure&\text{Count}&
 \pi_1&\text{Params}\cr
\exstrut&&&&&\cr
\noalign{\hrule}
\exstrut&&&&&\cr
1&\bD_{11}\splus\bA_4\splus2\bA_2&
 \ref{fig.d}(a)&
 (1,0)&
 \tref{s.RR0}&\cr
2&\bD_{9}\splus2\bA_4\splus\bA_2&
 \ref{fig.d}(a)&
 (1,0)&
 \tref{s.1.edge}&(5,5;3)\cr
3&\bD_{19}&
 \fragment(b)1&
 (1,0)&
 \tref{s.RR0}&\cr
4&\same1 \bD_{7}\splus\bA_{12}&
 \fragment(b){1}& 
 (1,0)&
 \tref{s.1.edge}&(13,\0;1)\cr
5&\same1 \bD_{7}\splus\bA_{12}&
 \fragment(b){1}&
 (0,1)&
 \tref{s.2.edge}&(13,13;1)\cr
6&\bD_{17}\splus\bA_2&
 \fragment(b){2}&
 (1,0)&
 \tref{s.RR0}&\cr
7&\bD_{9}\splus\bA_{10}&
 \fragment(b){2}&
 (1,0)&
 \tref{s.RR0}&\cr
8&\bD_{7}\splus\bA_{10}\splus\bA_2&
 \fragment(b)2&
 (0,1)&
 \tref{s.3.edge}&(3,11;1)\cr
9&\bD_{13}\splus\bA_6&
 \fragment(c){1,\bar1\!\!\!}&
 (0,1)&
 \tref{s.RR0}&\cr
10&\bD_{7}\splus2\bA_6&
 \fragment(c){1,\bar1\!\!\!}&
 (0,1)&
 \tref{s.2.edge}&(7,7;1)\cr
11&\bD_{15}\splus\bA_4&
 \fragment(c)2&
 (1,0)&
 \tref{s.RR0}&\cr
12&\bD_{11}\splus\bA_8&
 \fragment(c)2&
 (1,0)&
 \tref{s.1.edge}&(9,\0;5)\cr
13&\same2 \bD_{7}\splus\bA_8\splus\bA_4&
 \fragment(c)2&
 (1,0)&
 \tref{s.1.edge}&(9,\0;1)\cr
14&\same2 \bD_{7}\splus\bA_8\splus\bA_4&
 \fragment(c)2&
 (1,0)&
 \tref{s.3.edge}&\*(9,\0;1)\cr
15&\bD_{13}\splus\bA_4\splus\bA_2&
 \fragment(c)3&
 (1,0)&
 \tref{s.RR0}&\cr
16&\bD_{11}\splus\bA_6\splus\bA_2&
 \fragment(c)3&
 (1,0)&
 \tref{s.2.edge}&(3,7;5)\cr
17&\bD_{9}\splus\bA_6\splus\bA_4&
 \fragment(c)3&
 (1,0)&
 \tref{s.RR0}&\cr
18&\bD_{7}\splus\bA_6\splus\bA_4\splus\bA_2&
 \fragment(c)3&
 (1,0)&
 \tref{s.18}&(3,7;1)\cr
\exstrut&&&&&\cr
\noalign{\hrule}
\crcr}}}
\endinsert

\Remark\label{rem.regions}
Items~$4$ and~$5$ in Table~\ref{tab.d} differ by the choice of the
monogonal region~$\Rii$. We assume that \Nos$4$ and~$5$ correspond,
respectively, to the regions marked with~$\Ga$ or $\Gb$,
$\bar\Gb$ in Figure~\ref{fig.d}(b). (In the latter case, the two
choices differ by an orientation reversing symmetry, \ie, the two
curves are complex conjugate.)
Similarly, we assume that \Nos$13$ and~$14$ in the table
correspond,
respectively, to the monogonal regions marked with~$\Ga$ and~$\Gb$
in Figure~\ref{fig.d}(c).
\endRemark

\subsection{The case~$\bD_5$}\label{s.D5}
In this case, $\Sk$ has eight \black-vertices and two
distinguished regions, a bigon~$\Ri$ and a quadrilateral~$\Rii$.
If $\Ri$ is adjacent to~$\Rii$, then the two regions form together
an insertion shown in grey in Figure~\ref{fig.insertion}(b);
we call this fragment a \emph{bibigon}. As in the previous
subsection, removing the insertion and patching together the two
black edges, one obtains a $3$-regular skeleton~$\Sk'$ with four
\black-vertices. The new skeleton~$\Sk'$
represents an irreducible curve if and
only if so does~$\Sk$; hence, $\Sk'$ is one of the three
skeletons shown in Figure~\ref{fig.d}. Varying the position of the
insertion, one obtains items~$19$--$26$ in Table~\ref{tab.d5}.

\Remark
Unlike Subsection~\ref{s.D7}, this time the insertion has a
certain orientation, which should be taken into account. For this
reason, some positions shown in Figure~\ref{fig.d} give rise to
two rows in the table. Similarly, most positions shown in
Figure~\ref{fig.d-r} below give rise to two rows in
Table~\ref{tab.d5-r}.
\endRemark

\midinsert
\table\label{tab.d5}
Maximal sets of singularities with a type~$\bD_5$ point
\endtable

\def\no{}
\def\*{\llap{$^*$}}
\def\same{\afterassignment\dosame\count0=}
\def\dosame{\llap{$^{\the\count0\,}$}}
\def\tref#1{\,\text{\ref{#1}}\,}
\centerline{\vbox{\offinterlineskip\halign{%
\tabstrut\ \hss#\hss\ \vrule
&\quad$#$\hss\quad\vrule
&\quad\null#\hss\quad\vrule
&&\ \hss$#$\hss\ \vrule\cr
\noalign{\hrule}
\exstrut&&&&&\cr
\#&\text{Set of singularities}&
 \hss Figure&\text{Count}&
 \pi_1&\text{Params}\cr
\exstrut&&&&&\cr
\noalign{\hrule}
\exstrut&&&&&\cr
19&\bD_{5}\splus\bA_6\splus\bA_4\splus2\bA_2&
 \ref{fig.d}(a)&
 (1,0)&
 \tref{s.group.D5}&(7,5,3)\cr
20&\bD_{5}\splus\bA_{14}&
 \fragment(b)1&
 (0,1)&
 \tref{s.group.D5}&(15,\0,1)\cr
21&\bD_{5}\splus\bA_{12}\splus\bA_2&
 \fragment(b)2&
 (1,0)&
 \tref{s.group.D5}&(13,3,\0)\cr
22&\same3 \bD_{5}\splus\bA_{10}\splus\bA_4&
 \fragment(b)2&
 (1,0)&
 \tref{s.group.D5}&\,\*(11,5,\0)\,\cr
23&\bD_{5}\splus\bA_8\splus\bA_6&
 \fragment(c){1,\bar1\!\!\!}&
 (0,1)&
 \tref{s.group.D5}&\*(9,7,\0)\cr
24&\same3 \bD_{5}\splus\bA_{10}\splus\bA_4&
 \fragment(c){2}&
 (0,1)&
 \tref{s.group.D5}&(11,\0,1)\cr
25&\same4 \bD_{5}\splus\bA_8\splus\bA_4\splus\bA_2&
 \fragment(c){3}&
 (1,0)&
 \tref{s.group.D5}&(9,3,\0)\cr
26&\same5 \bD_{5}\splus\bA_6\splus2\bA_4&
 \fragment(c){3}&
 (1,0)&
 \tref{s.group.D5}&(5,7,\0)\cr
\exstrut&&&&&\cr
\noalign{\hrule}
\exstrut&&&&&\cr
27&\bD_{5}\splus(\bA_8\splus3\bA_2)&
 \frag1&
 (1,0)&
 \tref{s.group.27}&\cr
28&\bD_{5}\splus\bA_{10}\splus2\bA_2&
 \frag2&
 (1,0)&
 \tref{s.D5.1}&(11,\0,1)\cr
29&\same4 \bD_{5}\splus\bA_8\splus\bA_4\splus\bA_2&
 \frag{3,\bar3\!\!\!}&
 (0,1)&
 \tref{s.D5.1}&(5,9,3)\cr
30&\same5 \bD_{5}\splus\bA_6\splus2\bA_4&
 \frag4&
 (1,0)&
 \tref{s.D5.1}&(5,5,7)\cr
\exstrut&&&&&\cr
\noalign{\hrule}
\crcr}}}
\endinsert

Otherwise (if $\Ri$ is not adjacent to~$\Rii$), removing~$\Ri$
produces a skeleton~$\Sk'$ with six \black-vertices and a
distinguished quadrilateral region~$\Rii$. Such skeletons can
easily be classified; they are shown in Figure~\ref{fig.d5-r}
(where $\Rii$ is the outer region of the skeleton).
Using Theorem~\ref{th.splitting}, one can see that only one of
these skeletons
(the last one in Figure~\ref{fig.d5-r},
also shown in Figure~\ref{fig.d5}) represents an irreducible
curve. Varying the position of the bigonal insertion~$\Ri$ (shown
in grey and numbered in Figure~\ref{fig.d5}), one obtains
items~$27$--$30$ in Table~\ref{tab.d5}.

\midinsert
\centerline{\cpic{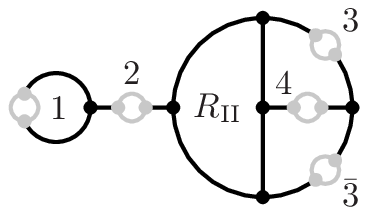}}
\figure\label{fig.d5}
A type~$\bD_5$ singular point: irreducible curves
\endfigure
\endinsert

\midinsert
\centerline{\cpic{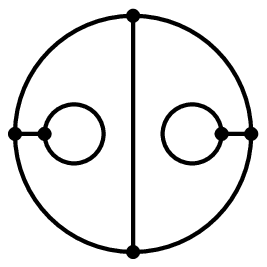}\qquad
\cpic{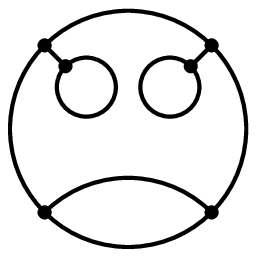}\qquad
\cpic{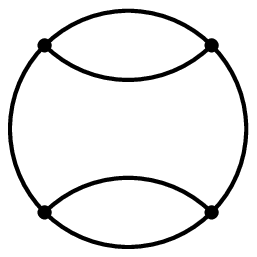}\qquad
\cpic{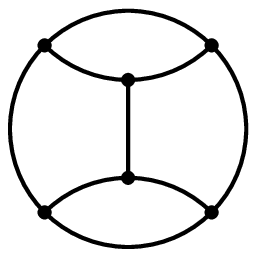}}
\bigskip
\centerline{\cpic{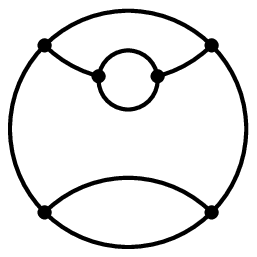}\qquad
\cpic{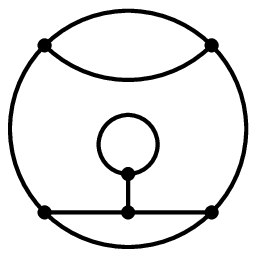}\qquad
\cpic{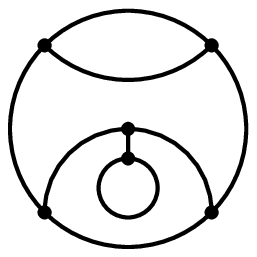}\qquad
\cpic{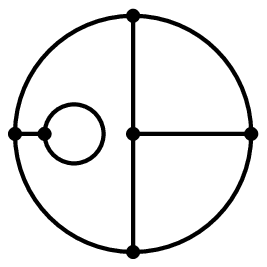}}
\figure\label{fig.d5-r}
A type~$\bD_5$ singular point: all curves
\endfigure
\endinsert

\Remark
The classification in this subsection could as well be
obtained from~\cite{Beukers}, where all $3$-regular skeletons with
eight \black-vertices are listed.
\endRemark

\section{Proof of Theorem~\ref{th.group}\label{S.proof.group}}

We fix a sextic $\CC\subset\Cp2$ as in the theorem and consider the
fundamental group
$\pi_1:=\pi_1(\Cp2\sminus\CC)=\pi_1(\Sigma_1\sminus(\BB\cup E))$,
where $\BB\subset\Sigma_1$ is the generalized trigonal curve
obtained by blowing up the type~$\bD$ singular point of~$\CC$,
\cf\ Subsection~\ref{s.models}. The group $\pi_1$ is computed on a
case by case basis, using the approach of Section~\ref{S.slopes}
and the skeletons found in
Section~\ref{S.proof.d}. (We retain the notation $\Ri$, $\Rii$ for
the two distinguished regions containing the improper fibers.)
Without further references, finite groups are treated using
\GAP~\cite{GAP}: in most cases, the {\tt Size} function
returns {\tt 6}, which suffices to conclude that the group
is~$\CG6$ (as so is its abelianization).

\subsection{A singular point of type~$\bD_p$, $p\ge7$}\label{s.group.D7}
We take for the reference fiber~$F_b$ the fiber over an
appropriate vertex~$\vi$ in the boundary of~$\Ri$ and choose a
canonical basis $\{\Ga_1,\Ga_2,\Ga_3\}$ in~$F$ corresponding to
the marking
at~$\vi$ such that
$\corner12\subset\Ri$.
Next, we choose an appropriate vertex~$\vii$ in the
boundary of~$\Rii$ and a canonical basis $\{\Gb_1,\Gb_2,\Gb_3\}$
in the fiber over~$\vii$ corresponding to a marking at~$\vii$
such that
$\corner23\subset\Rii$.
The translation
homomorphisms from~$\vi$ to~$\vii$
are computed below on a case by case basis.

According to~\iref{E.smooth}{smooth.12}
and~\iref{E.typeA}{typeA.23}, the slopes over~$\Ri$ and~$\Rii$ are
$\slope\ixi=\Ga_1\Ga_2$ and $\slope\ixii=\Gb_1$, respectively,
and the corresponding braid relations are
$$
[\Ga_1\Ga_2,\Ga_3]=1,\qquad
\Gs_2^{p-6}(\Gb_j)=\Gb_1\Gb_j\Gb_1\1,\quad j=2,3.
\eqtag\label{rel.D7}
$$
Furthermore, in view of the first relation in~\eqref{rel.D7},
the relation at infinity~\eqref{rel.infinity} simplifies to
$$
(\Ga_1\Ga_2)^2\Ga_3^3=\Gb_1.
\eqtag\label{inf.D7}
$$
These four relations are present in any group~$\pi_1$.

\subsection{The case of $\Rii$ adjacent to~$\Ri$}\label{s.RR0}
Assume that the region~$\Rii$ is adjacent to the insertion~$\Ri$
(\Nos$1$, $3$, $6$, $7$, $9$, $11$, $15$, and~$17$ in
Table~\ref{tab.d}).
Then $\vi$ and~$\vii$ can be chosen to coincide,
so that $\Gb_j=\Ga_j$, $j=1,2,3$,
and the relation at infinity~\eqref{inf.D7} simplifies
further to $\Ga_2\Ga_1\Ga_2\Ga_3^3=1$. It follows that
$\Ga_3$ commutes with $\Ga_2\Ga_1\Ga_2$; hence, in view of the
first relation in~\eqref{rel.D7}, one has
$[\Ga_3,\Ga_2]=[\Ga_3,\Ga_1]=1$. On the other hand,
$\Ga_1=\Ga_2\1\Ga_3^{-3}\Ga_2\1$ belongs to the abelian subgroup
generated by~$\Ga_2$ and~$\Ga_3$. Thus, the group is abelian.

The argument above applies to a reducible maximizing
sextic~$\CC$ as well,
provided that~$\CC$ is covered by Proposition~\ref{1-1.d} (\ie, $\CC$
has a type~$\bD_p$, $p\ge7$, singular point~$P$ and has no linear
components through~$P$) and the distinguished regions~$\Ri$,
$\Rii$ of the skeleton~$\Sk$ of the trigonal model~$\B$ of~$\CC$ are
adjacent to each other. Such curves can easily be enumerated
similar to Subsection~\ref{s.D7}, by reducing~$\Sk$ to a
skeleton~$\Sk'$
with at most four \black-vertices (see \eg~\cite{symmetric} and
Figure~\ref{fig.d-r};
note that this time we do not require that $\Sk'$ should have
exactly four \black-vertices as we accept curves~$\B$ with $\tD{}$
type singular fibers). The resulting sets of singularities are
listed in Table~\ref{tab.d-r}.

\midinsert
\centerline{\vbox{\halign{\hss#\hss&&\qquad\hss#\hss\cr
\cpic{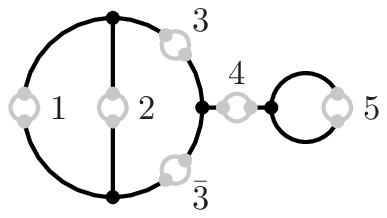}&\cpic{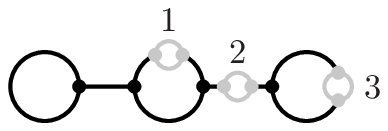}\cr
\noalign{\medskip}
(a)&(b)\cr
\crcr}}}
\bigskip
\centerline{\vbox{\halign{\hss#\hss&&\kern3em \hss#\hss\cr
\cpic{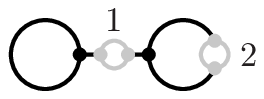}&\cpic{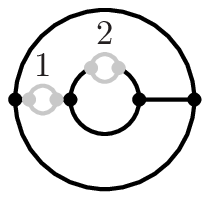}&\cpic{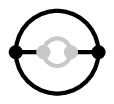}\cr
\noalign{\medskip}
(c)&(d)&(e)\cr
\crcr}}}
\figure\label{fig.d-r}
A type $\bD_p$ singular point, $p\ge7$: reducible curves
\endfigure
\endinsert

\midinsert
\def\fragment(#1){\ref{fig.d-r}(#1)\FRAG}
\table\label{tab.d-r}
Some reducible sextics with abelian fundamental groups
\endtable
\centerline{\vbox{\offinterlineskip\halign{%
\tabstrut\quad$#$\hss\quad\vrule
&\quad\null#\hss\quad\vrule\cr
\noalign{\hrule}
\exstrut&\cr
\text{Set of singularities}&\hss Figure\cr
\exstrut&\cr
\noalign{\hrule}
\exstrutii&\cr
\multispan2\tabstrut\qquad The splitting $C_3+C_3$\hss\vrule\cr
\exstrutii&\cr
\noalign{\hrule}
\exstrut&\cr
\bD_{14}\splus\bA_3\splus\bA_2&\fragment(a)1\cr 
\bD_{10}\splus\bA_7\splus\bA_2&\fragment(a)1\cr 
\bD_{16}\splus\bA_2\splus\bA_1&\fragment(a)4\cr 
\bD_{18}\splus\bA_1&\fragment(b)2\cr 
\bD_{14}\splus\bD_5&\fragment(c)1\cr 
\exstrut&\cr
\noalign{\hrule}
\exstrutii&\cr
\multispan2\tabstrut\qquad The splitting $C_4+C_2$\hss\vrule\cr
\exstrutii&\cr
\noalign{\hrule}
\exstrut&\cr
\bD_{11}\splus\bA_5\splus\bA_3&\fragment(a)2\cr 
\bD_{11}\splus\bA_7\splus\bA_1&\fragment(a)3\cr 
\bD_{9}\splus\bA_7\splus\bA_2\splus\bA_1&\fragment(a)5\cr 
\bD_{10}\splus\bA_9&\fragment(b)1\cr 
\bD_9\splus\bA_9\splus\bA_1&\fragment(b)3\cr 
\bD_{10}\splus\bD_9&\fragment(c)2\cr 
\bD_9\splus\bD_5\splus\bA_5&\fragment(c)2\cr 
\exstrut&\cr
\noalign{\hrule}
\crcr}}\,
\vbox{\offinterlineskip\halign{%
\tabstrut\quad$#$\hss\quad\vrule
&\quad\null#\hss\quad\vrule\cr
\noalign{\hrule}
\exstrut&\cr
\text{Set of singularities}&\hss Figure\cr
\exstrut&\cr
\noalign{\hrule}
\exstrutii&\cr
\multispan2\tabstrut\qquad The splitting $C_5+C_1$\hss\vrule\cr
\exstrutii&\cr
\noalign{\hrule}
\exstrut&\cr
\bD_{10}\splus\bA_5\splus\bA_4&\fragment(a)2\cr 
\bD_{14}\splus\bA_4\splus\bA_1&\fragment(a)3\cr 
\bD_{14}\splus2\bA_2\splus\bA_1&\fragment(a)5\cr 
\bD_{16}\splus\bA_3&\fragment(b)1\cr 
\bD_{16}\splus\bA_2\splus\bA_1&\fragment(b)3\cr 
\bD_{12}\splus\bD_7&\fragment(c)2\cr 
\bD_{12}\splus\bD_5\splus\bA_2&\fragment(c)2\cr 
\exstrut&\cr
\noalign{\hrule}
\exstrutii&\cr
\multispan2\tabstrut\qquad The splitting $C_3+C_2+C_1$\hss\vrule\cr
\exstrutii&\cr
\noalign{\hrule}
\exstrut&\cr
\bD_{12}\splus\bA_5\splus2\bA_1&\fragment(d)1\cr
\bD_{12}\splus2\bA_3\splus\bA_1&\fragment(d)2\cr
\bD_{10}\splus\bA_5\splus\bA_3\splus\bA_1&\fragment(d)2\cr
\bD_{10}\splus\bD_8\splus\bA_1&\ref{fig.d-r}(e)\cr
\bD_{10}\splus\bD_6\splus\bA_3&\ref{fig.d-r}(e)\cr
\exstrut&\cr
\noalign{\hrule}
\crcr}}}
\endinsert

\Remark
In Figures~\ref{fig.d-r}(c) and~(e), in addition to~$\Rii$ one
should also choose one of the remaining regions to contain the
$\tD{}$ type singular fiber (as one should have $\sum\tdfB=1$, see
Remark~\ref{rem.td}).
The degrees of the components of~$C$
can be determined using Theorem~\ref{th.splitting}, which
describes the components in terms of generators
(see Remark~\ref{rem.splitting}), and the
abelianization of relation~\eqref{inf.D7}.
\endRemark

\subsection{The case of $\Ri$ and~$\Rii$ connected by a single
edge}\label{s.1.edge}
Assume that $\Ri$ and~$\Rii$ are connected by a single edge of the
skeleton
(\Nos$2$, $4$, $12$, and~$13$ in Table~\ref{tab.d})
and choose vertices~$\vi$ and~$\vii$ so that this
connecting edge is
$[\vi,\vii]$. Then the translation homomorphism~$\gm_{[3,1]}$,
see~\ref{Sk.tr},
is given by
$$
 \Gb_1=\Ga_2\Ga_3\Ga_2\1,\quad
 \Gb_2=\Ga_2,\quad
 \Gb_3=(\Ga_2\Ga_3)\1\Ga_1(\Ga_2\Ga_3).
\eqtag\label{eq.1edge}
$$
Furthermore, in addition to~\eqref{rel.D7} and~\eqref{inf.D7}, one
has the braid relations
$$
\{\Ga_1,\Ga_3\}_l=\{\Ga_2,\Ga_3\}_m=1
\eqtag\label{braid.D7}
$$
from the two regions adjacent to~$\Ri$, see
Corollary~\ref{region.relations}.
(In most cases, these two
regions coincide. Since $[\Ga_1\Ga_2,\Ga_3]=1$, see~\eqref{rel.D7},
the relation $\{\Ga_1,\Ga_2\Ga_3\Ga_2\1\}_l=1$ given by the
corollary is equivalent to $\{\Ga_1,\Ga_3\}_l=1$.)
Trying the possible values of $(l,m;p)$, see
Table~\ref{tab.d},
with
\GAP~\cite{GAP}, one concludes that all four groups are abelian.

\subsection{The case of $\Ri$ and~$\Rii$ connected by two edges}\label{s.2.edge}
Assume that $\Ri$ and~$\Rii$ are connected by a chain~$\zeta$
of two edges
(\Nos$5$, $10$, and~$16$ in Table~\ref{tab.d}) and choose
reference vertices~$\vi$
and~$\vii$ at the two ends of~$\zeta$. The translation homomorphism
$\gm_{[3,1]}\circ\gm_{[3,1]}$
along~$\zeta$ is given by
$$
 \Gb_1=\Gd\1\Ga_1\Gd,\quad
 \Gb_2=\Ga_2,\quad
 \Gb_3=\Gr\1\Gd\Gr.
$$
where $\Gd=\Ga_2\Ga_3\Ga_2\1$, and relations~\eqref{rel.D7},
\eqref{inf.D7}, and~\eqref{braid.D7} with the values of $(l,m;p)$
given in Table~\ref{tab.d} suffice to show that all three groups
are abelian.

\subsection{The case of $\Ri$ and~$\Rii$ connected by three edges}\label{s.3.edge}
Assume that $\Ri$ and~$\Rii$ are connected by a chain~$\zeta$
of three edges
(\Nos$8$ and~$14$ in Table~\ref{tab.d})
and choose~$\vi$ and~$\vii$ at the ends of~$\zeta$.
Under an appropriate choice of~$\zeta$, the
translation homomorphism
$\gm_{[2,1]}\circ\gm_{[3,1]}\circ\gm_{[3,1]}$
is
$$
 \Gb_1=\Gd\1\Gr\Ga_2\Gr\1\Gd,\quad
 \Gb_2=\Gd\1\Ga_1\Gd,\quad
 \Gb_3=\Gr\1\Gd\Gr,
$$
where $\Gd=\Ga_2\Ga_3\Ga_2\1$. (To make this homomorphism uniform,
for \No$8$ we take for~$\Rii$ the monogon marked with~$\Gb$ in
Figure~\ref{fig.d}(b). Since the two curves in \No$8$ are complex
conjugate, their groups are isomorphic.)

For \No$8$, relations~\eqref{rel.D7},
\eqref{inf.D7}, and~\eqref{braid.D7} with $(l,m;p)=(3,11;1)$
suffice to show that the group is abelian.
For \No$14$, one should also take into account the relation
$\{\Gb_1,\Gb_2\}_5=1$ resulting from the pentagon adjacent
to~$\Rii$.

\subsection{The remaining case: \No18 in Table~\ref{tab.d}}\label{s.18}
In this case, the regions~$\Ri$ and $\Rii$ are connected by a
chain~$\zeta$ of
four edges. Choosing~$\vi$ and~$\vii$ at the ends of~$\zeta$, one
obtains the translation homomorphism
$\gm_{[2,1]}\circ\gm_{[3,1]}\circ\gm_{[3,1]}\circ\gm_{[3,1]}$
given by
$$
 \Gb_1=\Gd\1\Ga_1\1\Gd\Gr\Ga_2\Gr\1\Gd\1\Ga_1\Gd,\quad
 \Gb_2=\Gd\1\Ga_1\1\Gd\Ga_1\Gd,\quad
 \Gb_3=\Gr\1\Gd\1\Ga_1\Gd\Gr,
\eqtag\label{eq.18}
$$
where $\Gd=\Ga_2\Ga_3\Ga_2\1$.
Relations~\eqref{rel.D7},
\eqref{inf.D7}, and~\eqref{braid.D7} with $(l,m;p)=(3,7;1)$
suffice to show that the group is abelian.

\subsection{A singular point of type~$\bD_5$,
a bibigonal insertion}\label{s.group.D5}
In the case of a type~$\bD_5$ singular point, choose a pair of
reference
vertices~$\vi$ and~$\vii$ and canonical bases
$\{\Ga_1,\Ga_2,\Ga_3\}$ over~$\vi$ and $\{\Gb_1,\Gb_2,\Gb_3\}$
over~$\vii$
similar to
Subsection~\ref{s.group.D7}.

\midinsert
\centerline{\vbox{\halign{\hss#\hss&&\qquad\qquad\hss#\hss\cr
\cpic{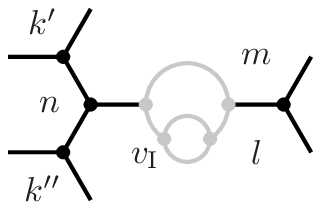}&
\cpic{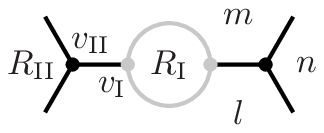}\cr
\noalign{\medskip}
(a)&(b)\cr
\crcr}}}
\figure\label{fig.regions}
Regions used in the computation
\endfigure
\endinsert

If the skeleton~$\Sk$ of~$\B$ has a bibigonal insertion, $\vi$
and~$\vii$ can be chosen to coincide, see
Figure~\ref{fig.regions}(a), so that
one has $\Gb_j=\Ga_j$, $j=1,2,3$.
According to~\ref{E.smooth}, the
slopes over~$\Ri$ and~$\Rii$ are
$\slope\ixi=\Ga_1\Ga_2$ and $\slope\ixii=(\Ga_2\Ga_3)^2$, respectively,
the braid
relations become
$$
[\Ga_1\Ga_2,\Ga_3]=[\Ga_1,(\Ga_2\Ga_3)^2]=1,
\eqtag\label{rel.D5}
$$
and the relation at infinity~\eqref{rel.infinity} simplifies to
$$
\Ga_3\Gr^3=(\Ga_2\Ga_3)^2.
\eqtag\label{inf.D5}
$$
Besides, $\pi_1$ has extra relations
$$
\{\Ga_1,\Ga_3\}_l=\{\Ga_2,\Gr\1\Ga_1\Gr\}_m=
\{\Ga_2,\Gr\1\Ga_2\Ga_3\Ga_2\1\Gr\}_n=1
\eqtag\label{rel.D5.ex}
$$
resulting from the $l$-, $m$-, and $n$-gonal regions marked in the
figure.

From the second relation in~\eqref{rel.D5} it follows that
$(\Ga_2\Ga_3)^2$ commutes with~$\rho$; then \eqref{inf.D5} implies
that $(\Ga_2\Ga_3)^2$ also commutes with~$\Ga_3$ and hence
with~$\Ga_2$. Thus, $(\Ga_2\Ga_3)^2$ is a central element and we
replace~$\pi_1$ with its quotient $G:=\pi_1/(\Ga_2\Ga_3)^2$.
(Otherwise,
the coset enumeration may fail in~\GAP.) Using \GAP~\cite{GAP}, we
show that $G=\CG2$; then $\pi_1$ is a central extension of a
cyclic group, hence abelian.

In most cases, for the conclusion that $G=\CG2$ it suffices to use
relations~\eqref{rel.D5}--\eqref{rel.D5.ex} with the values of
parameters $(l,m,n)$ listed in Table~\ref{tab.d5}.
For \Nos$22$ and~$23$ in Table~\ref{tab.d5},
one should also take into account the
relations
$$
\gather
\Ga_2=\Gr\1\Ga_2\Gr\1\Ga_2\Ga_3\Ga_2\1\Gr\Ga_2\1\Gr
 \quad\text{(for \No$22$)},\\
\Ga_2=\Gr\1\Ga_2\Gr\1\Ga_1\Gr\Ga_2\1\Gr
 \quad\text{(for \No$23$)}
\endgather
$$
resulting from appropriate monogonal regions of~$\Sk$.

This computation applies as well to a reducible maximizing
sextic~$\CC$, provided that it is covered by
Proposition~\ref{1-1.d5}, the skeleton~$\Sk$ of the trigonal
model~$\B$ of~$\CC$ has a bibigonal insertion, and $\CC$ splits into
two components (Figures~\ref{fig.d-r}(a)--(c)
and Table~\ref{tab.d5-r}; the latter
condition assures that the abelianization $G/[G,G]$ is finite).
This time, the central element $(\Ga_2\Ga_3)^2$ has infinite order
in the abelianization $\pi_1/[\pi_1,\pi_1]$ and hence one has
$[\pi_1,\pi_1]=[G,G]$.

\midinsert
\def\fragment(#1){\ref{fig.d-r}(#1)\FRAG}
\table\label{tab.d5-r}
Some reducible sextics with a type~$\bD_5$ point
\endtable
\def\*{\llap{$^*$}}
\centerline{\vbox{\offinterlineskip\halign{%
\tabstrut\quad$#$\hss\quad\vrule
 &\quad\null#\hss\quad\vrule&\ \hss$#$\hss\ \vrule
 &\quad$#$\ \ \hss\vrule\cr
\noalign{\hrule}
\exstrut&&&\cr
\text{Set of singularities}&\hss Figure&\text{Params}
 &[\pi_1,\pi_1]\cr
\exstrut&&&\cr
\noalign{\hrule}
\exstrutii\span\omit\span\omit&\cr
\multispan4\tabstrut\qquad The splitting $C_3+C_3$
 \ ($G/[G,G]=\CG4\oplus\CG3$)\hss\vrule\cr
\exstrutii\span\omit\span\omit&\cr
\noalign{\hrule}
\exstrut&&&\cr
\bD_5\splus\bA_7\splus\bA_5\splus\bA_2
&\fragment(a)1&
(6,8,3;\0,\0)
&\SL(2,\F_7)\cr
\bD_5\splus\bA_9\splus\bA_3\splus\bA_2
&\fragment(a)1&
(10,4,3;\0,\0)
&\{1\}\cr
\bD_5\splus\bA_{11}\splus\bA_2\splus\bA_1
&\fragment(a)4&
(12,12,1;\0,\0)
&\{1\}\cr
\bD_5\splus\bA_{13}\splus\bA_1
&\fragment(b)2&
(14,14,1;\0,\0)
&\{1\}\cr
\bD_{14}\splus\bD_5
&\fragment(c)1&
\*(\0,\0,1;\0,\0)
&\{1\}\cr
2\bD_5\splus\bA_9
&\fragment(c)1&
(10,10,1;\0,\0)
&\CG5\cr
\exstrut&&&\cr
\noalign{\hrule}
\exstrutii\span\omit\span\omit&\cr
\multispan4\tabstrut\qquad The splitting $C_4+C_2$
 \ ($G/[G,G]=\CG2\oplus\CG2$)\hss\vrule\cr
\exstrutii\span\omit\span\omit&\cr
\noalign{\hrule}
\exstrut&&&\cr
\bD_5\splus2\bA_5\splus\bA_4
&\fragment(a)2&
(6,5,6;\0,\0)
&\text{see~\ref{s.ex.2}}\cr
\bD_5\splus\bA_9\splus\bA_4\splus\bA_1
&\fragment(a){3,\bar3\!\!\!}&
(10,5,2;\0,\0)
&\{1\}\cr
\bD_5\splus\bA_9\splus2\bA_2\splus\bA_1
&\fragment(a)5&
(10,3,10;\0,2)
&\{1\}\cr
\bD_5\splus\bA_{11}\splus\bA_3
&\fragment(b)1&
(12,4,12;\0,1)
&\CG4\cr
\bD_5\splus\bA_{11}\splus\bA_2\splus\bA_1
&\fragment(b)3&
(12,3,12;\0,2)
&\text{see~\ref{s.ex.3}}\cr
\bD_{12}\splus\bD_5\splus\bA_2
&\fragment(c)2&
(\0,3,\0;\0,1)
&\{1\}\cr
\bD_7\splus\bD_5\splus\bA_7
&\fragment(c)2&
(8,\0,8;\0,1)
&\CG4\cr
2\bD_5\splus\bA_7\splus\bA_2
&\fragment(c)2&
(8,3,8;\0,\0)
&\text{see~\ref{s.ex.4}}\cr
\exstrut&&&\cr
\noalign{\hrule}
\exstrutii\span\omit\span\omit&\cr
\multispan4\tabstrut\qquad The splitting $C_5+C_1$
 \ ($G/[G,G]=\CG8$)\hss\vrule\cr
\exstrutii\span\omit\span\omit&\cr
\noalign{\hrule}
\exstrut&&&\cr
\bD_5\splus\bA_6\splus\bA_5\splus\bA_3
&\fragment(a)2&
(7,4,6;\0,\0)
&\{1\}\cr
\bD_5\splus\bA_7\splus\bA_6\splus\bA_1
&\fragment(a){3,\bar3\!\!\!}&
(7,8,2;\0,\0)
&\{1\}\cr
\bD_5\splus\bA_7\splus\bA_4\splus\bA_2\splus\bA_1\!\!
&\fragment(a)5&
(5,8,8;2,\0)
&\{1\}\cr
\bD_5\splus\bA_9\splus\bA_5
&\fragment(b)1&
(6,10,10;1,\0)
&\{1\}\cr
\bD_5\splus\bA_9\splus\bA_4\splus\bA_1
&\fragment(b)3&
(5,10,10;2,\0)
&\{1\}\cr
\bD_{10}\splus\bD_5\splus\bA_4
&\fragment(c)2&
(5,\0,\0;1,\0)
&\{1\}\cr
\bD_9\splus\bD_5\splus\bA_5
&\fragment(c)2&
(\0,6,6;1,\0)
&\{1\}\cr
2\bD_5\splus\bA_5\splus\bA_4
&\fragment(c)2&
(5,6,6;\0,\0)
&\{1\}\cr
\exstrut&&&\cr
\noalign{\hrule}
\crcr}}}
\endinsert

For some curves, we also take into account the additional
relations
$$
\{\Gr\1\Gd\Gr,\Gd\1\Ga_1\Gd\}_{k'}=1
\quad\text{or}\quad
\{\Gr\Ga_2\Gr\1,\Gd\1\Ga_1\Gd\}_{k''}=1,
\eqtag\label{rel.D5.ex-r}
$$
$\Gd=\Ga_2\Ga_3\Ga_2\1$,
resulting from the $k'$- and $k''$-gonal regions marked in
Figure~\ref{fig.regions}(a).
As usual, we always skip the relations corresponding to a region
of~$\Sk$ containing a type~$\tD{}$ singular fiber of~$\B$ (as
these relations differ from those indicated above).

For most curves,
relations~\eqref{rel.D5}--\eqref{rel.D5.ex-r} with
the values of $(l,m,n;k',k'')$ given in Table~\ref{tab.d5-r}
suffice to identify~$\pi_1$, either because the group is already
abelian or due to Proposition~\ref{prop.-1}. The few special cases
are discussed below.

\subsubsection{The set of singularities
$\bD_{14}\splus\bD_5$}\label{s.ex.1}
To show that $\pi_1$ is abelian, one needs to take into account
the additional relation
$$
(\Ga_2\Ga_3)\Ga_2(\Ga_2\Ga_3)\1=\Gr\Gd\Gr\1
$$
resulting from the other monogonal region of~$\Sk$.

\subsubsection{The set of singularities
$\bD_5\splus2\bA_5\splus\bA_4$}\label{s.ex.2}
In this case, \GAP~\cite{GAP} shows that
the commutant $[\pi_1,\pi_1]$ is one of the five
perfect groups of order~$7680$. I do not know which of the five
groups it is.

\subsubsection{The set of singularities
$\bD_5\splus\bA_{11}\splus\bA_2\splus\bA_1$}\label{s.ex.3}
One has $[\pi_1,\pi_1]=\Z$. Although
$[\pi_1,\pi_1]$
infinite,
it can be simplified using the \GAP\ commands
\vskip\abovedisplayskip
\halign{\qquad\tt#\hss\cr
P := PresentationNormalClosure(g, Subgroup(g, [g.1/g.2]));\cr
SimplifyPresentation(P);\cr}
\vskip\belowdisplayskip\noindent
which return a presentation with a single generator and no
relations.

\subsubsection{The set of singularities
$2\bD_5\splus\bA_7\splus\bA_2$}\label{s.ex.4}
One has
$$
\pi_1'/\pi_1''=\Z\oplus\CG3
\quad\text{and}\quad
\pi_1''=Q_8:=\{\pm1,\pm i,\pm j,\pm k\}\subset\Bbb H,
$$
where we abbreviate $\pi_1'=[\pi_1,\pi_1]$ and
$\pi_1''=[\pi_1',\pi_1']$. The first statement is straightforward.
For the second one, consider the normal closure~$H$
of~$\Ga_1\Ga_2\1$ in~$G$. It is an order~$3$ subgroup of~$G'$ and one
has
$H/[H,H]=\Z$.
Hence $[H,H]=[G',G']=G''$.
Simplifying the presentation of~$H$ given by~\GAP\
in the same way as in~\ref{s.ex.3}, one
obtains two generators $\Gk_2$, $\Gk_3$ and three relations
$$
\Gk_3^4=1,\quad
\Gk_3^2\Gk_2\1\Gk_3^2\Gk_2=1,\quad
\Gk_2^{-2}\Gk_3\1\Gk_2\Gk_3\1\Gk_2\Gk_3=1.
$$
From the first two relations it follows that $\Gk_3^2$ is a
central element. Then the third one, rewritten in the form
$\Gk_2^3=(\Gk_2\Gk_3\1)^3\Gk_3^2$, implies that $\Gk_2^3$ is also
central. Since the image of~$\Gk_2^3$ in $H/[H,H]$ has infinite
order, $[H,H]$ is equal to the commutant of the quotient
$H/\Gk_2^3$. The latter group is finite
and its commutant is~$Q_8$.

\proposition\label{perturbations}
Let~$\CC$ be one of the sextics listed in Tables~\ref{tab.d-r}
and~\ref{tab.d5-r},
and let~$\CC'$ be an irreducible perturbation of~$\CC$ preserving
the distinguished type~$\bD_5$ singular point. Then
$\pi_1(\Cp2\sminus\CC')=\CG6$.
\endproposition

\proof
It suffices to consider one of the seven sextics~$\CC$
in Table~\ref{tab.d5-r} that have
nonabelian groups. Since the distinguished type~$\bD_5$ point is
preserved, the perturbation of~$\CC$ is followed by a perturbation
$\B\to\B'$ of its trigonal model, or a perturbation $\BB\to\BB'$
of the generalized trigonal curve in~$\Sigma_1=\Cp2(P)$
which is actually used in the computation of~$\pi_1$:
one can assume that a proper
type~$\bA_{2r-1}$ or~$\bD_{2r-1}$ singular point~$Q$ of~$\BB$ is
perturbed so that the intersection $\BB'\cap\MB_Q$ is connected,
where $\MB_Q$ is a Milnor ball about~$Q$.

If the point~$Q$ that is perturbed is of type~$\bD_{2r-1}$,
the inclusion homomorphism
$\pi_1(\MB_Q\sminus\BB)\to\pi_1(\Sigma_1\sminus\BB)$
is onto (as $\MB_Q$ contains all three generators in a fiber
sufficiently close to~$Q$).
On the other hand, for any
perturbation $\BB\to\BB'$ with $\BB'\cap\MB_Q$ connected, the
group $\pi_1(\MB_Q\sminus\BB')$ is abelian (see~\cite{dessin.e8};
the maximal perturbation with the connectedness property is
$\bD_{2r-1}\to\bA_{2r-2}$.)

If $Q$ is of type~$\bA_{2r-1}$, the group of~$\BB'$ is found
similar to that of~$\BB$: it suffices to replace the corresponding
(necessarily even) parameter(s) in $(l,m,n;k',k'')$, see
Table~\ref{tab.d5-r}, with its maximal odd divisor. Considering
curves and parameters one by one and using \GAP~\cite{GAP}, one
concludes that all groups are abelian.
\endproof

\subsection{Other curves not of torus type}\label{s.D5.1}
In all three cases (\Nos$28$--$30$ in Table~\ref{tab.d5}),
the distinguished regions~$\Ri$ and~$\Rii$ are connected by
a single edge, see Figure~\ref{fig.d5}.
Choose vertices~$\vi$, $\vii$ introduced in
Subsection~\ref{s.group.D5}
as shown in
Figure~\ref{fig.regions}(b). Then the translation
homomorphism from~$\vi$ to~$\vii$ is given by~\eqref{eq.1edge}.
Hence, the braid relations from~$\Ri$ and~$\Rii$ and the relation
at infinity become
$$
[\Ga_1\Ga_2,\Ga_3]=[(\Gr\1\Ga_1\Gr\Ga_2)^2,\Ga_3]=1
$$
and
$$
\Ga_3\Gr^3=(\Ga_2\Gr\1\Ga_1\Gr)^2,
$$
respectively. Consider also the relations
$$
\{\Ga_1,\Ga_3\}_l=\{\Ga_2,\Ga_3\}_m=
\{\Ga_2\1\Ga_1\Ga_2,\Gr\1\Ga_2\Gr\}_n=1
$$
resulting from the $l$-, $m$-, and $n$-gonal regions marked in
Figure~\ref{fig.regions}(b). Using the values of $(l,m,n)$ given
in Table~\ref{tab.d5}, one concludes that all three groups are
abelian.

\subsection{The curve of torus type
(\No27 in Table~\ref{tab.d5})}\label{s.group.27}
The skeleton~$\Sk$ is the one shown in Figure~\ref{fig.d5}, with
the bigonal insertion labelled with~$1$. Take for~$\vi$
and~$\vii$, respectively, the upper vertex of the insertion and
the center of the large circle in the figure. Then a complete set
of relations for~$\pi_1$ is
$$
\gathered
[\Ga_1\Ga_2,\Ga_3]=\{\Ga_2,\Ga_3\}_3=\{\Ga_1,\Ga_3\}_9=1,\\
[(\Gb_1\Gb_2)^2,\Gb_3]=\{\Gb_2,\Gb_3\}=\{\Gb_1,\Gb_2\Gb_3\Gb_2\1\}_3=1,\\
\Ga_3\Gr^3=(\Gb_1\Gb_2)^2,
\endgathered
$$
where $\Gb_1$, $\Gb_2$, $\Gb_3$ are as in~\eqref{eq.18}.
(We use the marking at~$\vii$
such that $\corner12\subset\Rii$. Along an appropriate path of
length~$4$, the translation homomorphism from~$\vi$ to~$\vii$ is
$\gm_{[2,1]}\circ\gm_{[3,1]}\circ\gm_{[3,1]}\circ\gm_{[3,1]}$.)
Using the \GAP\ commands
\vskip\abovedisplayskip
\halign{\qquad\tt#\hss\cr
P := PresentationNormalClosure(g, Subgroup(g, [g.1/g.3]));\cr
SimplifyPresentation(P);\cr}
\vskip\belowdisplayskip\noindent
one finds that $[\pi_1,\pi_1]$ is a free group on two generators.
Since there is a canonical (perturbation) epimorphism
$\pi_1\onto\bar\BG3$ and all groups involved are residually
finite, hence Hopfian, the above epimorphism is an isomorphism.
(This approach to using \GAP~\cite{GAP} to treat a group
`suspected' to be isomorphic to~$\bar\BG3$
was suggested to me by E.~Artal Bartolo.)


\widestnumber\key{EO1}
\refstyle{C}


\widestnumber\no{99}
\Refs

\ref{AVG}
\by V.~I.~Arnol$'$d, A.~N.~Varchenko, S.~M.~Guse\u{\i}n-Zade
\book Singularities of differentiable maps
\vol I. The classification of critical points, caustics and wave fronts
\publ Nauka
\publaddr Moscow
\yr     1982
\lang Russian
\transl\nofrills English translation:
\book Monographs in Mathematics
\vol 82
\publ Birkh\"auser Boston, Inc.
\publaddr Boston, MA
\yr 1985
\endref\label{AVG}

\ref{A}
\by E.~Artal Bartolo
\paper A curve of degree five with non-abelian fundamental group
\jour Topology Appl.
\vol 79
\yr 1997
\issue 1
\pages 13--29
\endref\label{Artal}

\ref{BM}
\by F.~Beukers, H.~Montanus
\paper Explicit calculation of elliptic fibrations of $K3$-surfaces
and their Belyi-maps
\inbook Number theory and polynomials
\pages 33--51
\bookinfo London Math. Soc. Lecture Note Ser.
\vol 352
\publ Cambridge Univ. Press
\publaddr Cambridge
\yr 2008
\endref\label{Beukers}

\ref{CM}
\by H.~S.~M.~Coxeter, W.~O.~J.~Moser
\book  Generators and relations for discrete groups
\yr  1957
\publ Springer--Verlag
\publaddr Berlin-G\"ottingen-Heidelberg
\endref\label{Coxeter}

\ref{D3}
\by A.~Degtyarev
\paper Quintics in $\C\roman{p}^2$ with nonabelian fundamental group
\jour Algebra i Analis
\yr 1999
\vol    11
\issue  5
\pages  130--151
\lang Russian
\moreref\nofrills English transl. in
\jour Leningrad Math.~J.
\vol 11
\yr 2000
\issue 5
\pages 809--826
\endref\label{groups}

\ref{D8}
\by A.~Degtyarev
\paper Fundamental groups of symmetric sextics
\jour  J. Math. Kyoto Univ.
\vol 48
\issue 4
\yr 2008
\pages 765--792
\endref\label{degt.e6}

\ref{D7}
\by A.~Degtyarev
\paper Stable symmetries of plane sextics
\jour Geometri{\ae} Dedicata
\vol 137
\yr 2008
\issue 1
\pages 199--218
\endref\label{symmetric}

\ref{D5}
\by A.~Degtyarev
\paper Zariski $k$-plets via dessins d'enfants
\jour Comment. Math. Helv.
\vol 84
\issue 3
\yr 2009
\pages 639--671
\endref\label{degt.kplets}

\ref{D7}
\by A.~Degtyarev
\paper Plane sextics via dessins d'enfants
\finalinfo\tt arXiv:0812.3258
\endref\label{dessin.e7}

\ref{D6}
\by A.~Degtyarev
\paper Plane sextics with a type $\bE_8$ singular point
\finalinfo\tt arXiv:0902.2281
\endref\label{dessin.e8}

\ref{D6}
\by A.~Degtyarev
\paper Plane sextics with a type $\bE_6$ singular point
\finalinfo\tt arXiv:0907.4714
\endref\label{dessin.e6}

\ref{D7}
\by A.~Degtyarev
\paper Topology of plane algebraic curves\rom: the algebraic approach
\finalinfo\tt arXiv:0907.0289
\endref\label{survey}

\ref{D7}
\by A.~Degtyarev
\paper Transcendental lattice of an extremal elliptic surface
\finalinfo\tt arXiv:0907.1809
\endref\label{tripods}

\ref{GAP}
\by The GAP Group
\book GAP---Groups, Algorithms, and Programming
\bookinfo Version 4.4.10
\yr 2007
\finalinfo ({\tt http:\allowbreak//www.gap-system.org})
\endref\label{GAP}

\ref{vK}
\by E.~R.~van~Kampen
\paper On the fundamental group of an algebraic curve
\jour  Amer. J. Math.
\vol   55
\yr    1933
\pages 255--260
\endref\label{vanKampen}

\ref{Ko}
\by K.~Kodaira
\paper On compact analytic surfaces, II--III
\jour Annals of Math.
vol 77--78
\yr 1963
\pages 563--626, 1--40
\endref\label{Kodaira}

\ref{OP1}
\by M.~Oka, D.~T.~Pho
\paper Classification of sextics of torus type
\jour Tokyo J. Math.
\vol 25
\issue 2
\pages 399--433
\yr 2002
\endref\label{OkaPho.moduli}

\ref{Sh}
\by I.~Shimada
\paper On the connected components of the moduli of
polarized $K3$ surfaces
\toappear
\endref\label{Shimada}

\ref{Ya}
\by J.-G.~Yang
\paper Sextic curves with simple singularities
\jour Tohoku Math. J. (2)
\vol 48
\issue 2
\yr 1996
\pages 203--227
\endref\label{Yang}

\ref{Z1}
\by O.~Zariski
\paper On the problem of existence of algebraic functions of two
variables possessing a given branch curve
\jour Amer. J. Math.
\vol 51
\yr 1929
\pages 305--328
\endref\label{Zariski.group}

\endRefs
\enddocument